# NON-ASYMPTOTIC THEORY FOR NONPARAMETRIC TESTING

By Yun Yang[*], Zuofeng Shang[†] and Guang Cheng[‡]

*Florida State University, Binghamton University and Purdue University*
February 3, 2017

We consider nonparametric testing in a non-asymptotic framework. Our statistical guarantees are exact in the sense that Type I and II errors are controlled for any finite sample size. Meanwhile, one proposed test is shown to achieve minimax optimality in the asymptotic sense. An important consequence of this non-asymptotic theory is a new and practically useful formula for selecting the optimal smoothing parameter in nonparametric testing. The leading example in this paper is smoothing spline models under Gaussian errors. The results obtained therein can be further generalized to the kernel ridge regression framework under possibly non-Gaussian errors. Simulations demonstrate that our proposed test improves over the conventional asymptotic test when sample size is small to moderate.

**1. Introduction.** Many classical statistical inferential procedures are built upon large sample theory that relies on a growing amount of data information or a large number of samples. However, in practice, it is often the case that only a small to moderate amount of samples are available, which limits the applicability of the classic asymptotic inferential procedures. Recently, finite sample inferential theory has provided a new perspective for statistical analysis. For instance, with the aid of various notions of concentration inequality, [2, 3, 17, 22, 23, 30] have developed statistical inference procedures that are theoretically valid for any fixed sample size. As far as we are aware, the parameter of interest in these works is finite dimensional. The goal of this paper is to develop finite sample theory in another important setting – nonparametric models.

In the asymptotic regime, there is a vast amount of literature devoted to developing theories for nonparametric inferences such as testing and confidence band; see [26, 28]. However, little progress has been gained towards finite sample theory for nonparametric inference procedures. Our work can be viewed as an initial attempt to establishing a non-asymptotic framework for nonparametric testing, which covers the existing asymptotic theory as a direct consequence. This effort requires new technical tools such as (uniformly valid) large deviation bounds. In particular, two Wald-type test statistic are constructed with their cut-off values being adjusted according to any finite sample. This is in sharp contrast with conventional asymptotic tests that rely on the null limit distributions, e.g., likelihood ratio test in [8, 20].

---

[*]Assistant Professor.

[†]Assistant Professor. Research Sponsored by a startup grant from Binghamton University.

[‡]Corresponding Author. Professor. Research Sponsored by NSF CAREER Award DMS-1151692, DMS-1418042, and Office of Naval Research.

We begin with smoothing spline models under Gaussian errors. As a first attempt, we consider a test statistic simply defined as a discrepancy between the null function and smoothing spline estimate, and develop a corresponding nonasymptotic deviation inequality (see Theorem 3.1). Despite its simplicity, this naive test is shown minimax sub-optimal in an asymptotic sense according to [12, 15, 4]. This is due to a non-vanishing bias term that reduces the power. This motivates the development of a more refined non-asymptotic deviation bound (see Theorem 3.4) for a "de-biased" smoothing spline estimate, based on which our second test statistic is proposed and shown to achieve the minimx optimality by correcting the bias explicitly. To our knowledge, this deviation bound is the first non-asymptotic expansion of the smoothing spline estimate up to a higher-order decaying remainder term. Based on these large deviation results, Type I and II errors are controlled for any finite sample size in both methods; see Corollaries 3.3 and 3.6. As an interesting by-product, a lower bound on the sample size is quantified to yield desirable Type I and II errors. The relation between the proposed tests and the asymptotically valid likelihood ratio test ([20]) is also highlighted. As for other smoothing-spline-based tests such as LMP ([7]), they were developed from Bayesian setup and always rely on complicated null distributions involving nuisance parameters, as reviewed in [16].

Three non-trivial generalizations are further carried out. First, we consider a general class of composite hypothesis, e.g., the null set is a subspace of polynomial functions. Second, we extend our test to more general regression settings where errors might be non-Gaussian. Third, we demonstrate that our finite sample testing framework also accommodates the more general kernel ridge regression (KRR) [21], including smoothing spline as a special case.

In practice, the choice of smoothing parameter is crucially important. However, it is known that the widely used generalized cross validation (GCV, [25]) does not lead to a minimax optimal testing procedure; see [8]. Importantly, our non-asymptotic theory yields a practically useful formula to select an optimal smoothing parameter. To be more specific, it is obtained by directly minimizing a separation function (that quantifies the minimal detectable signal strength) derived in the non-asymptotic framework. As far as we are aware, this non-asymptotic formula is new.

Our simulations demonstrate the empirical advantage of the above selection method over GCV in both the proposed test and the conventional asymptotic test such as penalized likelihood ratio test (PLRT, [20]). Additionally, the simulation study confirms that our finite-sample based testing procedure is uniformly more powerful than PLRT as the sample size grows (under the same choice of smoothing parameter). Besides the non-asymptotic design, another reason for this empirical success is that the proposed test removes a bias term from the PLRT due to penalization. Although this bias term is asymptotically of higher order, it can significantly deteriorate the power of the test when sample size is small to moderate. We count this as another highlight by applying our finite sample framework.

The rest of this paper is organized as follows. In Section 2, background and preliminaries are introduced. Sections 3 and 4 include the main results of this paper. In Section 3, test statistic based on smoothing spline estimation are constructed and their non-asymptotic properties are investigated. The methods are valid for both simple and composite hypotheses. Section 4 includes extensions to more general regression settings. Section 5 includes a simulation study. Main proofs are provided in Section 6, while the rest are deferred to a supplementary document [29].

## 2. Preliminary.
In this section, we state the nonparametric models and hypothesis of interest, and also review some basic theory of reproducing kernel Hilbert space (RKHS).

Let $\mathbf{D}_n = \{(Y_i, X_i) : i = 1, 2, \ldots, n\}$ be $iid$ random samples following the nonparametric regression model

$$(2.1) \qquad Y = f(X) + \epsilon,$$

where $\epsilon$ is a random error with mean zero and unit variance. For simplicity, we consider the one-dimensional case where $X \in \mathbb{I} := [0, 1]$, and assume that the design $X$ and $\epsilon$ are independent. Let the marginal density function of $X$ be $\pi(x)$ which satisfies $0 < \inf_{x \in \mathbb{I}} \pi(x) \le \sup_{x \in \mathbb{I}} \pi(x) < \infty$. We assume that $f$ belongs to an $m$-th order Sobolev space

$$S^m(\mathbb{I}) = \left\{ f \in L^2(\mathbb{I}) \middle| f^{(j)} \text{ are abs. cont. for } j = 0, 1, \ldots, m-1, \text{ and } \int_0^1 |f^{(m)}(x)|^2 dx < \infty \right\}.$$

Let $P_f$ denote the probability measure under $f$, and $E_f$ be the corresponding expectation.

Consider the following hypothesis testing problem:

$$(2.2) \qquad H_0 : f \in \mathcal{F}_0, \ vs \ H_1 : f \notin \mathcal{F}_0,$$

where $\mathcal{F}_0$ is a proper subset of $S^m(\mathbb{I})$. For example, $\mathcal{F}_0 = \{f_0\}$ for some known $f_0 \in S^m(\mathbb{I})$ in simple hypothesis, while $\mathcal{F}_0 = \{\text{all linear functions in } S^m(\mathbb{I})\}$ in composite hypothesis. Our testing statistic is constructed based on the following smoothing spline estimator

$$\widehat{f}_{n,\lambda} = \operatorname{argmax}_{f \in S^m(\mathbb{I})} \ \ell_{n,\lambda}(f),$$

where $\ell_{n,\lambda}(f)$ is the penalized loss function defined as

$$\ell_{n,\lambda}(f) = -\frac{1}{2n} \sum_{i=1}^n (Y_i - f(X_i))^2 - \frac{\lambda}{2} J(f, f),$$

where $\lambda > 0$ is the penalty parameter and $J(f, g) = \int_0^1 f^{(m)}(x) g^{(m)}(x) dx$ for any $f, g \in S^m(\mathbb{I})$.

Let $V(f, g) = E\{f(X)g(X)\}$ and $\langle f, g \rangle = V(f, g) + \lambda J(f, g)$. It follows that $\langle \cdot, \cdot \rangle$ defines a valid inner product in $S^m(\mathbb{I})$; see [20]. Endowed with $\langle \cdot, \cdot \rangle$, $S^m(\mathbb{I})$ is an RKHS. We use the notation $\| \cdot \|$ to denote the corresponding RKHS norm. Let $K(x_1, x_2)$ denote the reproducing kernel function, a function from $\mathbb{I} \times \mathbb{I}$ to $\mathbb{R}$ satisfying the reproducing property $\langle K_x, f \rangle = f(x)$ for any $x \in \mathbb{I}$ and $f \in S^m(\mathbb{I})$, where $K_x(\cdot) := K(x, \cdot)$ is an element in $S^m(\mathbb{I})$ for any $x \in \mathbb{I}$. Let $\mathcal{P}_\lambda$ denote a self-adjoint operator from $S^m(\mathbb{I})$ to itself that satisfies $\langle f, \mathcal{P}_\lambda g \rangle = \lambda J(f, g)$, for all $f, g \in S^m(\mathbb{I})$; see [20] for the existence of such an operator.

We make the following assumption on the existence of eigen-pairs $(\rho_\nu, \varphi_\nu)$ that simultaneously diagonalize $V$ and $J$. This assumption is commonly made in smoothing spline literature.

ASSUMPTION A1. *There exists $\varphi_\nu \in S^m(\mathbb{I})$, for $\nu \geq 1$, satisfying $c_\varphi := \sup_{\nu \geq 1} \|\varphi_\nu\|_{\sup} < \infty$, where $\|\cdot\|_{\sup}$ denotes the supremum norm, and a nondecreasing sequence of eigenvalues $\rho_\nu \sim (c_\rho \nu)^{2m1}$, where $c_\rho > 0$ is a constant, such that*

$$(2.3) \qquad V(\varphi_\mu, \varphi_\nu) = \delta_{\mu\nu}, \quad J(\varphi_\mu, \varphi_\nu) = \rho_\mu \delta_{\mu\nu}, \ \mu, \nu = 1, 2, \ldots,$$

*where $\delta_{\mu\nu}$ is the Kronecker's delta. Furthermore, any $g \in S^m(\mathbb{I})$ admits a Fourier expansion $g = \sum_\nu V(g, \varphi_\nu)\varphi_\nu$ with convergence held in the $\|\cdot\|$-norm.*

It follows from [20] that for any $x \in \mathbb{I}$, $K_x = \sum_{\nu \geq 1} \frac{\varphi_\nu(x)}{1+\lambda\rho_\nu}\varphi_\nu$. Let $h = \lambda^{1/(2m)}$. Under Assumption A1, it is straightforward to verify the following property of $K_x$

$$(2.4) \qquad \|K_x\|^2 = K(x,x) = \sum_{\nu \geq 1} \frac{\varphi_\nu(x)^2}{1+\lambda\rho_\nu} \lesssim h^{-1}.$$

## 3. Main Results.
In this section, we construct two nonparametric test methods based on $\widehat{f}_{n,\lambda}$. The first type is straightforward but sub-optimal in the minimax sense, while the second attains the minimax optimality by removing the bias from the former. Our major contribution is to develop non-asymptotic theoretical analysis for both testing methods. Specifically, Type I and II errors can be controlled for any *finite* sample size. This leads to a non-asymptotic formula in selecting the optimal $\lambda$ in practice. These non-asymptotic results are developed based on large deviation bounds between $\widehat{f}_{n,\lambda}$ and $f$ that are shown uniformly valid over an "unit ball" in $S^m(\mathbb{I})$.

In this section, we assume Gaussian errors and uniform design, i.e., $\epsilon \sim N(0,1)$ and $X \sim \mathrm{Unif}(\mathbb{I})$, and postpone extensions to general error distributions to the next section.

### 3.1. *A Preliminary First-Order Testing Procedure.*
To illustrate the idea, we first focus on the simple hypothesis $H_0 : f = f_0$, where $f_0$ is a known function in $S^m(\mathbb{I})$. Under the null, $\widehat{f}_{n,\lambda}$ converges to $f_0$ under suitable norms as sample size $n$ tends to infinity. Naturally, the deviation between $\widehat{f}_{n,\lambda}$ and $f_0$ can be used as a test statistic:

$$(3.1) \qquad T_{n,\lambda} = \|\widehat{f}_{n,\lambda} - f_0\|.$$

Despite the simple form of $T_{n,\lambda}$, it will be shown to be asymptotically sub-optimal. We call $T_{n,\lambda}$ as first-order testing in the rest of the paper.

The null hypothesis is rejected if $T_{n,\lambda}$ exceeds some threshold $d_n(M)$ (to be defined later), where $M$ is a constant controlling the significance level of the test. To calibrate a finite sample valid $d_n(M)$, we first need to establish a large deviation bound for $T_{n,\lambda}$ uniformly over an "unit ball" in $S^m(\mathbb{I})$, defined as $H^m(1) = \{f \in S^m(\mathbb{I}) | J(f,f) \leq 1\}$.

Let $c_K = \sup_{x \in \mathbb{I}} \sqrt{hK(x,x)}$, which is finite by following (2.4). Recall that $h = \lambda^{1/(2m)}$.

---

[1]For two positive sequences $a_\nu$ and $b_\nu$, denote $a_\nu \lesssim b_\nu$ if $a_\nu = O(b_\nu)$; $a_\nu \gtrsim b_\nu$ if $b_\nu = O(a_\nu)$; $a_\nu \asymp b_\nu$ if $a_\nu \lesssim b_\nu$ and $b_\nu \lesssim a_\nu$; $a_\nu \sim b_\nu$ if $a_\nu/b_\nu$ tends to one when $\nu \to \infty$.

THEOREM 3.1. *Let Assumption A1 be satisfied. For any positive constants $(h, r, M)$ satisfying $c_K^2 \sqrt{M} r h^{-1/2} A(h) \leq 1/2$, it holds that*

$$(3.2) \qquad \sup_{f \in H^m(1)} P_f\left(\|\widehat{f}_{n,\lambda} - f\| \geq \delta_n(M, r)\right) \leq 2\exp(-Mnhr^2),$$

*where $\delta_n(M, r) = 2h^m + c_K(\sqrt{2M}r + (nh)^{-1/2})$ and $A(h)$ is an explicit function of $h$ defined in Remark 3.1 below.*

Theorem 3.1 is proven in Section 6.2.

REMARK 3.1. *The function $A(h)$ in Theorem 3.1 is defined as $A(h) = A(h, 2)$. As will be seen below, $A(h, \varepsilon)$ is an explicit formula depending on Dudley's entropy integral, which controls the upper bound of a concentration inequality (Lemma 6.1) that plays a key role in the proof of Theorem 3.1. Specifically,*

$$
\begin{aligned}
A(h, \varepsilon) \;=\; & \frac{32\sqrt{6}}{\tau} c_K^{-1} c_0^m h^{-(2m-1)/2} \Psi\left(\frac{1}{2} c_K c_0^{-m} h^{(2m-1)/2} \varepsilon\right) \\
& + \frac{20\sqrt{6}\varepsilon}{\tau} \sqrt{\log\left(1 + \exp\left(2c_0(c_K h^{(2m-1)/2}\varepsilon)^{-1/m}\right)\right)}, \quad \text{for } \varepsilon > 0,
\end{aligned}
$$

*where $\tau = \sqrt{\log 1.5} \approx 0.6368$ and the function $\Psi(r)$, resulted from Dudley's entropy integral (see [24]), is defined as $\Psi(r) = \int_0^r \sqrt{\log(1 + \exp(x^{-1/m}))}dx$. In the above, $c_0$ is chosen as the constant controlling the packing number of $\mathcal{G} := \{g \in S^m(\mathbb{I}) : \|g\|_{\sup} \leq 1, J(g, g) \leq c_K^{-2} h^{-2m+1}\}$, i.e., $c_0$ satisfies, for any $\varepsilon > 0$,*

$$(3.3) \qquad \log N(\varepsilon, \mathcal{G}, \|\cdot\|_{\sup}) \leq c_0 h^{-\frac{2m-1}{2m}} \varepsilon^{-1/m},$$

*where $N(\varepsilon, \mathcal{G}, \|\cdot\|_{\sup})$ is the $\varepsilon$-packing number. Existence of such $c_0$ follows from [24, 20].*

According to (3.2), we choose $r = (nh)^{-1/2}$. In this case, the threshold becomes $\delta_n(M, (nh)^{-1/2})$, denoted as $d_n(M)$. The following theorem, as an immediate consequence of Theorem 3.1, characterizes the upper bounds of Type I and II errors for any finite sample size.

THEOREM 3.2. *Let the Assumptions in Theorem 3.1 be satisfied, and let $M > 0$ and $L > 1$ be any constants. Given that $d_n(M) := \delta_n(M, (nh)^{-1/2}) = 2h^m + c_K(\sqrt{2M} + 1)(nh)^{-1/2}$, then it holds that*

$$
\begin{aligned}
\text{Type I error}: & \quad P_{f_0}(T_{n,\lambda} \geq d_n(M)) \leq 2\exp(-M), \\
\text{Type II error}: & \quad \sup_{\substack{f - f_0 \in H^m(1) \\ \|f - f_0\| \geq \rho_n(M, L)}} P_f(T_{n,\lambda} \leq d_n(M)) \leq 2\exp(-L),
\end{aligned}
$$

*where $\rho_n(M, L) = 4h^m + c_K(\sqrt{2M} + \sqrt{2L} + 2)(nh)^{-1/2}$.*

Theorem 3.2 implies that the Type I error falls below $\alpha$ and Type II error falls below $\beta$ if we set $M = M_0 := \log(2/\alpha)$ and $L = L_0 := \log(2/\beta)$, respectively. Consequently, the separation rate between $f$ and $f_0$ (in terms of $\|\cdot\|$), under which the testing power is at least $1 - \beta$, is $\rho_n(M_0, L_0)$.

We summarize the above discussions in the Corollary below.

COROLLARY 3.3. *Under the assumptions of Theorem 3.1, for any $\alpha, \beta \in (0, 1)$, we have*

$$\textit{Type I error}: \quad P_{f_0}(T_{n,\lambda} \geq d_n(\log(2/\alpha))) \leq \alpha,$$

$$\textit{Type II error}: \quad \sup_{\substack{f - f_0 \in H^m(1) \\ \|f - f_0\| \geq \rho_n(\log(2/\alpha), \log(2/\beta))}} P_f(T_{n,\lambda} \leq d_n(\log(2/\alpha))) \leq \beta,$$

*where $\rho_n(\log(2/\alpha), \log(2/\beta)) = 4h^m + c_K(\sqrt{2\log(2/\alpha)} + \sqrt{2\log(2/\beta)} + 2)(nh)^{-1/2}$, as a function of $h$, achieves its minimum at*

$$(3.4) \qquad h_* = \left( \frac{c_K^2(\sqrt{\log(2/\alpha)} + \sqrt{\log(2/\beta)} + \sqrt{2})^2}{8m^2 n} \right)^{1/(2m+1)}.$$

The value of $c_K = \sup_{x \in \mathbb{I}} \sqrt{hK(x, x)}$ in (3.4) can be approximately determined by (2.4) which in turn requires the estimate of $\|\varphi_\nu\|_{\sup}$. The latter estimation can be done numerically by the spectral decomposition of the reproducing kernel matrix w.r.t. $J$. For instance, the R packages *gss* ([9]) and *assist* ([27]) both allow us to extract the kernel matrix corresponding to $J$, and the eigenvectors of the matrix provide a good estimate of the eigenfunctions.

The minimal separation rate (computed at $h_*$) is given as

$$\rho_n(M_0, L_0) = D(c_K, \alpha, \beta) n^{-\frac{m}{2m+1}},$$

where $D(c_K, \alpha, \beta)$ is a positive constant depending on $c_K$ and $(\alpha, \beta)$ only. Nevertheless, the above rate fails to match with the minimax lower bound, namely, $n^{-4m/(4m+1)}$ (see [12, 15, 4]). Hence, the first order testing procedure is sub-optimal from an asymptotic perspective. A finite sample valid and asymptotically optimal testing method is further proposed in Section 3.2.

3.2. *An Optimal Second-Order Testing Procedure.* In this subsection, we improve the first order testing to attain minimax optimality. A closer examination of $T_{n,\lambda}$ reveals that its sub-optimality is due to a large bias arising from the deviation $\widehat{f}_{n,\lambda} - f_0$, which inflates the separation gap $\rho_n(M, L)$. Fortunately, this bias can be easily removed as shown in the following second order deviation result. This observation motivates a new testing procedure, i.e., (3.6).

THEOREM 3.4. *Let Assumption A1 be satisfied. For any positive constants $(h, r, M)$ satisfying $c_K^2 \sqrt{M} rh^{-1/2} A(h) \leq 1/2$, it holds that*

$$(3.5) \quad \sup_{f \in H^m(1)} P_f\left( \left\| \widehat{f}_{n,\lambda} - f - \left( \frac{1}{n} \sum_{i=1}^n \epsilon_i K_{X_i} - \mathcal{P}_\lambda f \right) \right\| \geq \gamma_n(M, r) \right) \leq 2\exp(-Mnhr^2),$$

*where $\gamma_n(M, r) = c_K^2 \sqrt{M} rh^{-1/2} A(h)\delta_n(M, r)$ and $\delta_n(M, r)$ is defined in Theorem 3.1.*

Theorem 3.4 is proven in Section 6.2 by employing a concentration inequality (see Lemma 6.1) for an operator-valued empirical process and a contraction mapping argument.

We remark that the threshold $\gamma_n(M, r)$ in Theorem 3.4 converges to zero faster than the threshold $\delta_n(M, r)$ in Theorem 3.1 if we choose $(r, h)$ to satisfy $r^2 = (nh)^{-1}$, $h = o(1)$ and $n^{-1}h^{-(6m-1)/(2m)} = o(1)$ (this leads to $rh^{-1/2}A(h) = o(1)$ by noting that $A(h) \lesssim h^{-(2m-1)/(4m)}$). In this case, $\|\widehat{f}_{n,\lambda} - f\|$ and $\|n^{-1}\sum_{i=1}^{n}\epsilon_i K_{X_i} - \mathcal{P}_\lambda f\|$ are of the same asymptotic order. For this reason, we call (3.2) a first-order deviation bound, and (3.5) a second-order deviation bound.

In view of (3.5), a second-order test statistic is developed as

$$(3.6) \qquad \widetilde{T}_{n,\lambda} = \|\widehat{f}_{n,\lambda} - (I - \mathcal{P}_\lambda)f_0\|^2 - \frac{1}{n}\sum_{\nu \geq 1}\frac{1}{1 + \lambda\rho_\nu},$$

where the second term $n^{-1}\sum_{\nu \geq 1}(1 + \lambda\rho_\nu)^{-1}$ is the expectation of $\|n^{-1}\sum_{i=1}^{n}\epsilon_i K_{X_i}\|^2$. The term $\mathcal{P}_\lambda f_0$ in $\widetilde{T}_{n,\lambda}$ is a bias correction term due to penalization. Subtracting a $f$-independent constant $n^{-1}\sum_{\nu \geq 1}(1 + \lambda\rho_\nu)^{-1}$ in $\widetilde{T}_{n,\lambda}$ is merely for technical simplicity in the subsequent derivations.

The corresponding testing rule is $\phi_{n,\lambda} = I(|\widetilde{T}_{n,\lambda}| \geq d_n(M, h))$, where $d_n(M, h)$ controls Type I error through the choice of $(M, h)$. Based on Theorem 3.4, we will prove in Theorem 3.5 that

$$(3.7) \qquad \begin{aligned} \text{Type I error}: \qquad & E_{f_0}\{\phi_{n,\lambda}\} \leq e_I(M), \\ \text{Type II error}: \qquad & \sup_{\substack{f - f_0 \in H^m(1) \\ \|f - f_0\| \geq \rho_n(M, L, h)}} E_f\{1 - \phi_{n,\lambda}\} \leq e_{II}(L), \end{aligned}$$

where $\rho_n(M, L, h)$, $e_I(M)$ and $e_{II}(L)$ are given in Theorem 3.5. Note that $d_n(M, h)$ and $\rho_n(M, L, h)$ are different from $d_n(M)$ and $\rho_n(M, L)$ defined in the previous section.

Under $H_0 : f = f_0$, we can decompose the test statistic $\widetilde{T}_{n,\lambda}$ as

$$\widetilde{T}_{n,\lambda} = \left\|\frac{1}{n}\sum_{i=1}^{n}\epsilon_i K_{X_i}\right\|^2 - \frac{1}{n}\sum_{\nu \geq 1}\frac{1}{1 + \lambda\rho_\nu} + \text{higher-order remainder}$$

$$(3.8) \qquad = \left[\frac{1}{n^2}\sum_{i,j=1}^{n}\epsilon_i\epsilon_j K(X_i, X_j) - E\left\{\frac{1}{n^2}\sum_{i,j=1}^{n}\epsilon_i\epsilon_j K(X_i, X_j)\right\}\right] + \text{higher-order remainder}.$$

By controlling the first two terms in (3.8), we can obtain large deviation bounds for $\widetilde{T}_{n,\lambda}$ under both null and alternative hypotheses (see Lemma 6.2 and Lemma 6.3 in Section 6.3). This leads to the following theorem characterizing the finite sample property of the proposed test $\widetilde{T}_{n,\lambda}$.

THEOREM 3.5. *Suppose Assumption A1 holds. For any constants $(h, M, L)$ satisfying $c_K^2\sqrt{M}$ $n^{-1/2}h^{-1}A(h) \leq 1/2$ and $c_K^2\sqrt{L}n^{-1/2}h^{-1}A(h) \leq 1/2$, we choose the cutoff value $d_n(M, h)$ as*

$$(3.9) \qquad d_n(M, h) = \frac{4\rho_K}{n\sqrt{h}}\sqrt{M} + R_{1,n}(M),$$

*where $\rho_K^2 = hE[K^2(X_1, X_2)]$ with $X_1, X_2 \overset{iid}{\sim} X$, and separation function*

$$(3.10) \qquad \rho_n(M, L, h) = \sqrt{\zeta_K\lambda} + \sqrt{\frac{2L}{n}} + \sqrt{d_n(M, h) + \frac{2L}{n} + R_{2,n}(L)},$$

*where $\zeta_K = \sup_{g \in H^m(1)} \lambda^{-1} \|\mathcal{P}_\lambda g\|^2$. Then (3.7) holds with*

$$
(3.11) \qquad\qquad e_I(M) = 15 \exp(-M) \qquad and \qquad e_{II}(L) = 30 \exp(-L).
$$

*Here explicit forms of the remainder terms $R_{1,n}(M)$ and $R_{2,n}(L)$ are provided in Section S.2 in the supplement.*

The following corollary is obtained as an immediate consequence of Theorem 3.5.

COROLLARY 3.6. *Under the assumptions of Theorem 3.4, for any $\alpha, \beta \in (0,1)$, we have*

*Type I error* : $\qquad\qquad P_{f_0}(\widetilde{T}_{n,\lambda} \geq d_n(\log(15/\alpha), h)) \leq \alpha,$

*Type II error* : $\qquad \sup_{\substack{f - f_0 \in H^m(1) \\ \|f - f_0\| \geq \rho_n(\log(15/\alpha), \log(30/\beta), h)}} P_f(\widetilde{T}_{n,\lambda} \leq d_n(\log(15/\alpha), h)) \leq \beta.$

An important implication of Theorem 3.5 is that $\widetilde{T}_{n,\lambda}$ is asymptotically minimax optimal. In fact, using the expression (3.9) of $d_n(M, h)$, we see that under the asymptotic regime $n \to \infty$, the leading term in $\rho_n(M, L, h)$ scales as

$$
\sqrt{\zeta_K \lambda} + \sqrt{\frac{4\rho_K}{n\sqrt{h}}\sqrt{M}}.
$$

By minimizing this leading term w.r.t. $h$, we obtain the minimal separation rate

$$
\rho_n(M, L, h_{**}) \asymp n^{-2m/(4m+1)}
$$

when $h$ is chosen as

$$
(3.12) \qquad\qquad h = h_{**} \equiv \left( \left(\frac{4\rho_K}{\zeta_K}\right)^2 M \right)^{1/(4m+1)} n^{-2/(4m+1)} \asymp n^{-2/(4m+1)}.
$$

The above minimal separation rate matches with the minimax rate of testing obtained in [12, 15, 4].

The practical implementation of $\widetilde{T}_{n,\lambda}$ requires us to estimate $(I - \mathcal{P}_\lambda)f_0$. Instead of direct estimation, we approximate it by the following "noiseless" version of smoothing spline estimator

$$
\widehat{f}_{n,\lambda}^{NL} = \operatorname{argmin}_{f \in S^m(\mathbb{I})} \frac{1}{n}\sum_{i=1}^{n}(f_0(X_i) - f(X_i))^2 + \lambda J(f, f).
$$

By applying Theorem 3.4 with $\epsilon_i \sim N(0,0)$, it is easy to see that $\|\widehat{f}_{n,\lambda}^{NL} - (I - \mathcal{P}_\lambda)f_0\|$ has the same derivation bound as $\left\|\widehat{f}_{n,\lambda} - f_0 - \left(\frac{1}{n}\sum_{i=1}^{n}\epsilon_i K_{X_i} - \mathcal{P}_\lambda f_0\right)\right\|$, and therefore

$$
P_{f_0}\left( \left\| \widehat{f}_{n,\lambda} - \widehat{f}_{n,\lambda}^{NL} - \frac{1}{n}\sum_{i=1}^{n}\epsilon_i K_{X_i} \right\| \geq 2\gamma_n(M, r) \right) \leq 4\exp(-Mnhr^2).
$$

As a consequence, we can replace $(I - \mathcal{P}_\lambda)f_0$ in $\widetilde{T}_{n,\lambda}$ by $\widehat{f}_{n,\lambda}^{NL}$. In simulations, this approximated version of $\widetilde{T}_{n,\lambda}$ (by using $\widehat{f}_{n,\lambda}^{NL}$) is found to work very well, and have a larger power than $T_{n,\lambda}$ especially when $f_0$ under the null is far from zero (so that $\mathcal{P}_\lambda f_0$ incurs a relatively large bias).

REMARK 3.2. *(Relation with likelihood ratio test) The non-asymptotic results obtained for $T_{n,\lambda}$ and $\widetilde{T}_{n,\lambda}$ can be extended to another type of nonparametric testing: likelihood ratio test. We first define the penalized likelihood ratio test (PLRT) as follows*

$$
\begin{aligned}
(3.13) \quad 2PLRT(g) \quad &:= \quad 2(\ell_{n,\lambda}(g) - \ell_{n,\lambda}(\widehat{f}_{n\lambda})) \\
&= \quad \frac{1}{n} \sum_{i=1}^{n} (\widehat{f}_{n,\lambda}(X_i) - g(X_i))^2 + \langle \widehat{f}_{n,\lambda} - g, \mathcal{P}_\lambda(\widehat{f}_{n,\lambda} - g) \rangle.
\end{aligned}
$$

*In comparison, our test statistic (up to constants) can be expressed as*

$$
(3.14) \qquad \|\widehat{f}_{n,\lambda} - g\|^2 = V(\widehat{f}_{n,\lambda} - g, \widehat{f}_{n,\lambda} - g) + \langle \widehat{f}_{n,\lambda} - g, \mathcal{P}_\lambda(\widehat{f}_{n,\lambda} - g) \rangle,
$$

*where $g = f_0$ for $T_{n,\lambda}$ and $g = (I - \mathcal{P}_\lambda)f_0$ for $\widetilde{T}_{n,\lambda}$. Note that the first term in (3.14) is the expectation of that in (3.13) according to the definition of $V(f, f)$. In Appendix S.6, we prove that for any $g \in S^m(\mathbb{I})$, the deviation between $2PLRT(g)$ and $\|\widehat{f}_{n,\lambda} - g\|^2$ is of higher order comparing to the dominating term of $d_n(M, h)$, i.e., $4\rho_K\sqrt{M}/(n\sqrt{h})$. Therefore, after some modifications, the results for $T_{n,\lambda}$ and $\widetilde{T}_{n,\lambda}$ also hold for $2PLRT(f_0)$ and $2PLRT((I - \mathcal{P}_\lambda)f_0)$, respectively.*

REMARK 3.3. *(Effective sample size for testing) An interesting consequence of Theorem 3.5 is that the minimal sample size, named as "effective" sample size, can be derived to simultaneously achieve the pre-determined size and power. We demonstrate this idea using a simple example. Suppose we want to test $H_0 : f = 0$ vs $H_1 : f = f_*$, where $f_* \in H^m(1)\backslash\{0\}$ and known. For any $\alpha, \beta \in (0, 1)$, choose $M = \log(15/\alpha)$ and $L = \log(30/\beta)$, and let $\rho_n(M, L, \hat{h}(M, L))$ be determined by (3.10). Here, $\hat{h}(M, L)$ is defined to minimize $\rho_n(M, L, h)$ for any fixed $M, L$. Theorem 3.5 says that the Type I and II errors have upper bounds $\alpha$ and $\beta$ respectively, if*

$$
(3.15) \qquad \|f_*\| \geq \rho_n(M, L, \hat{h}(M, L)).
$$

*Directly solving inequality (3.15) for the effective sample size is nontrivial. Instead we can get an approximate solution. Specifically, $\rho_n(M, L, \hat{h}(M, L))$ can be approximated by its leading term as*

$$
\begin{aligned}
&\rho_n(M, L, \hat{h}(M, L)) \\
&\approx \quad \zeta_K^{1/2} \left( \frac{4\zeta_K \rho_K \sqrt{M}}{n} \right)^{2m/(4m+1)} + \sqrt{\frac{2L}{n}} + \sqrt{\zeta_K \left( \frac{4\zeta_K \rho_K \sqrt{M}}{n} \right)^{4m/(4m+1)} + \frac{2L}{n}}.
\end{aligned}
$$

*Therefore, instead of solving (3.15), one can solve*

$$
(3.16) \quad \|f_*\| \geq \zeta_K^{1/2} \left( \frac{4\zeta_K \rho_K \sqrt{M}}{n} \right)^{2m/(4m+1)} + \sqrt{\frac{2L}{n}} + \sqrt{\zeta_K \left( \frac{4\zeta_K \rho_K \sqrt{M}}{n} \right)^{4m/(4m+1)} + \frac{2L}{n}}.
$$

*A numerical solution to (3.16) can be used as an "effective" sample size.*

3.3. *Extension to Composite Hypothesis.* Our results can be generalized to testing composite hypothesis. For example, one composite hypothesis of particular interest is whether $f$ is a polynomial function, say with degree less than $m$. In this case, a new test statistic is proposed with similar non-asymptotic guarantee and asymptotic minimax optimality. For technical simplicity, we assume $m \geq 2$ throughout this section.

For simplicity, we consider testing whether $f$ is linear in this section:

$$H_0 : f \in \mathcal{F}_0 \text{ vs } H_1 : f \notin \mathcal{F}_0,$$

where $\mathcal{F}_0 = \{f : f \text{ is linear on } \mathbb{I}\}$. We propose a test statistic as

$$\widetilde{T}_{n,\lambda}^{com} = \|\widehat{f}_{n,\lambda} - \widehat{f}_n^{H_0}\|^2 - \frac{1}{n} \sum_{\nu \geq 1} \frac{1}{1 + \lambda \rho_\nu},$$

where $\widehat{f}_n^{H_0}(x) = (1, x)(D_X^T D_X)^{-1} D_X^T Y$ is the maximum likelihood estimate under $H_0$. Here, $D_X = (\mathbf{1}, \mathbf{X})$ denotes the design matrix with intercept, and $\mathbf{X} = (X_1, \ldots, X_n)^T$. Note that $\mathcal{P}_\lambda \widehat{f}_n^{H_0} = 0$.

We propose to reject $H_0$ if and only if $|\widetilde{T}_{n,\lambda}^{com}| \geq d_n^{com}(M, h)$ for some threshold $d_n^{com}(M, h)$. In general, we will establish that for any finite sample size (as in (3.7))

$$
\begin{align}
(3.17) \quad & \sup_{f \text{ is linear}} P_f\left(|\widetilde{T}_{n,\lambda}^{com}| \geq d_n^{com}(M, h)\right) \leq \quad e_I^{com}(M), \\
& \sup_{\substack{f \in H^m(1) \\ \|f^\perp\| \geq \rho_n^{com}(M, L, h)}} P_f\left(|\widetilde{T}_{n,\lambda}^{com}| < d_n^{com}(M, h)\right) \leq \quad e_{II}^{com}(L),
\end{align}
$$

where $d_n^{com}(M, h), e_I^{com}(M), e_{II}^{com}(L)$ and $\rho_n^{com}(M, L, h)$ are given in Theorem 3.7. Here, $f^\perp$ denotes the projection of $f$ onto the orthogonal complement of the space of linear functions in $S^m(\mathbb{I})$, and satisfies $V(f^\perp, g) = E_X\{f^\perp(X)g(X)\} = 0$ for any linear function $g$ (hence, $E_X\{f^\perp(X)\} = E_X\{Xf^\perp(X)\} = 0$).

We still start from large deviation bounds of the test statistic under both null and alternative hypotheses. To do so, $\widetilde{T}_{n,\lambda}^{com}$ needs to be decomposed as (3.19) based on the following arguments. Let $f_0$ be the true linear function in $\mathcal{F}_0$ from which dataset $\mathbf{D}_n$ is generated, and let $f_\epsilon = \widehat{f}_n^{H_0} - f_0$. It follows by Theorem 3.4 that, under $H_0$, $\widetilde{T}_{n,\lambda}^{com}$ can be decomposed as

$$
\begin{align}
\widetilde{T}_{n,\lambda}^{com} &= \left\|\frac{1}{n} \sum_{i=1}^n \epsilon_i K_{X_i} - \mathcal{P}_\lambda f_0 - f_\epsilon\right\|^2 - \frac{1}{n} \sum_{\nu \geq 1} \frac{1}{1 + \lambda \rho_\nu} + \text{higher-order remainder} \\
&= \left[\frac{1}{n^2} \sum_{i,j=1}^n \epsilon_i \epsilon_j K(X_i, X_j) - E\left\{\frac{1}{n^2} \sum_{i,j=1}^n \epsilon_i \epsilon_j K(X_i, X_j)\right\}\right]
\end{align}
$$

$$(3.18) \quad + \|\mathcal{P}_\lambda f_0 + f_\epsilon\|^2 + \frac{2}{n} \sum_{i=1}^n \epsilon_i \mathcal{P}_\lambda f_0(X_i) + \frac{2}{n} \boldsymbol{\epsilon}^T D_X (D_X^T D_X)^{-1} D_X^T \boldsymbol{\epsilon} + \text{higher-order remainder},$$

where $\boldsymbol{\epsilon} = (\epsilon_1, \ldots, \epsilon_n)^T$. The second equality follows from the fact that

$$f_\epsilon(x) = (1, x)(D_X^T D_X)^{-1} D_X^T \boldsymbol{\epsilon} := \widehat{\alpha} + \widehat{\beta} x.$$

Since both $f_0$ and $f_\epsilon$ are linear functions, $\mathcal{P}_\lambda f_0 = \mathcal{P}_\lambda f_\epsilon = 0$, and hence, $\|f_\epsilon\|^2 = \|f_\epsilon\|_2^2 = \int_0^1 \left(f_\epsilon(x)\right)^2 dx$. As a consequence, the preceding display (3.18) can be further simplified as

$$
\begin{aligned}
\widetilde{T}_{n,\lambda}^{com} &= \Big[\frac{1}{n^2}\sum_{i,j=1}^n \epsilon_i \epsilon_j K(X_i, X_j) - E\Big\{\frac{1}{n^2}\sum_{i,j=1}^n \epsilon_i \epsilon_j K(X_i, X_j)\Big\}\Big] \\
&\quad + \|f_\epsilon\|_2^2 + \frac{2}{n}\boldsymbol{\epsilon}^T D_X (D_X^T D_X)^{-1} D_X^T \boldsymbol{\epsilon} + \text{higher-order remainder.}
\end{aligned}
\tag{3.19}
$$

Based on (3.19), we can control the type I and type II error of the test (see Lemma 6.4 and Lemma 6.5), yielding finite sample analysis for the composite hypothesis testing procedure.

THEOREM 3.7. *Suppose that Assumption A1 holds. For any constants* $(h, M, L)$ *satisfying* $c_K^2 \sqrt{M} n^{-1/2} h^{-1} A(h) \leq 1/2$ *and* $c_K^2 \sqrt{L} n^{-1/2} h^{-1} A(h) \leq 1/2$, *we choose*

$$
\begin{aligned}
d_n^{com}(M, h) &= \frac{4\rho_K}{n\sqrt{h}}\sqrt{M} + R_{1,n}^c(M), \quad \text{and} \\
\rho_n^{com}(M, L, h) &= \sqrt{\zeta_K \lambda} + \sqrt{\frac{2L}{n}} + \sqrt{d_n(M, h) + \frac{2L}{n} + R_{2,n}^c(L)}.
\end{aligned}
\tag{3.20}
$$

*Then (3.17) holds with* $e_I^{com}(M) = 24\exp(-M)$ *and* $e_{II}^{com}(L) = 60\exp(-L)$. *Here the forms of* $R_{1,n}^c(M)$ *and* $R_{2,n}^c(L)$ *are provided in Section S.2 in the supplement.*

Similarly, we have the following corollary for the composite test.

COROLLARY 3.8. *Under the assumptions of Theorem 3.7, for any* $\alpha, \beta \in (0, 1)$, *we have*

$$
\begin{aligned}
\textit{Type I error}: \quad &\sup_{f \text{ is linear}} P_f(\widetilde{T}_{n,\lambda}^{com} \geq d_n^{com}(\log(24/\alpha), h)) \leq \alpha, \\
\textit{Type II error}: \quad &\sup_{\substack{f \in H^m(1) \\ \|f^\perp\| \geq \rho_n^{com}(\log(24/\alpha),\log(60/\beta),h)}} P_f(\widetilde{T}_{n,\lambda}^{com} \leq d_n^{com}(\log(24/\alpha), h)) \leq \beta.
\end{aligned}
$$

Similar to the simple hypothesis testing, Theorem 3.7 provides a non-asymptotic approach to select the smoothing parameter $h$ by numerically minimizing the separation function $h \mapsto \rho_n^{com}(\log(24/\alpha), \log(60/\beta), h)$. Note that the leading terms of $d_n^{com}(M, h)$ and $\rho_n^{com}(M, L, h)$ are exactly the same as $d_n(M, h)$ and $\rho_n(M, L, h)$ in Section 3.2 for simple hypothesis. Hence, the selected $\hat{h}$ in both cases are asymptotically rate-equivalent. Under such an $\hat{h}$, the composite testing procedure in consideration is minimax rate-optimal.

## 4. Extensions.
In this section, we generalize our main results to non-Gaussian errors and the more general setup of kernel ridge regression.

### 4.1. General non-Gaussian Regression.
Suppose that $\mathbf{D}_n = \{(Y_i, X_i) : i = 1, 2, \ldots, n\}$ are *iid* samples generated from model (2.1), with errors $\epsilon_i$ whose log-likelihood is $\ell(\cdot)$. Suppose that the function $\ell(\cdot)$ is three-times continuously differentiable and is strictly concave. Let $\eta = \ell^{(1)}(\epsilon)$ and $\eta_i = \ell^{(1)}(\epsilon_i)$ for $i = 1, \ldots, n$, and $\sigma^2 = -E\{\ell^{(2)}(\epsilon)\} > 0$.

ASSUMPTION A2.  *There are positive constants $C_0, C_1$ such that*

$$(4.1) \qquad E\{\exp(|\ell^{(1)}(\epsilon)|/C_0)\} \le C_1, \ \ and \ \ E\left\{\exp\left(\sup_{|a|\le 1}|\ell^{(j)}(\epsilon+a)|/C_0\right)\right\} \le C_1, \ j = 2, 3.$$

*Furthermore, $E\{\eta\} = 0$ and $E\{\eta^2\} = \sigma^2$.*

Condition (4.1) says that $\ell^{(j)}(\epsilon)$ for $j = 1, 2, 3$ satisfies exponential tail condition. Note that for any $a \in \mathbb{R}$, $\int \exp(\ell(\epsilon + a))d\epsilon = 1$. Taking first- and second-order derivatives with respect to $a$, it can be shown that

$$(4.2) \qquad \int \exp(\ell(\epsilon+a))\ell^{(1)}(\epsilon+a)d\epsilon = \int \exp(\ell(\epsilon+a))[\ell^{(1)}(\epsilon+a)^2 + \ell^{(2)}(\epsilon+a)]d\epsilon = 0.$$

Setting $a = 0$ in (4.2), one gets $E\{\eta\} = 0$ and $E\{\eta^2\} = \sigma^2$. Therefore, Assumption A2 is a reasonable one.

Similar to Theorems 3.1 and 3.4 for the Gaussian errors, we develop the following derivation bounds for the first- and second-order approximations of the penalized likelihood estimate $\widehat{f}_{n,\lambda}$.

THEOREM 4.1.  *Suppose that Assumptions A1 and A2 are satisfied. For any positive $M, r_n, h$ satisfying Condition H in Section S.1, the following two results hold:*

*(i)*
$$\sup_{f \in H^m(1)} P_f\left(\|\widehat{f}_{n,\lambda} - f\| \ge \delta'_n(M)\right) \le (2C_1 + 4)\exp(-Mnhr_n^2),$$

*where $\delta'_n(M) = 2h^m + 24C_0 c_K(4C_1 + M)r_n$.*

*(ii)*
$$\sup_{f \in H^m(1)} P_f\left(\left\|\widehat{f}_{n,\lambda} - f - \frac{1}{n}\sum_{i=1}^{n}\eta_i K_{X_i} + \mathcal{P}_\lambda f\right\| > c_n(M)\right) \le (2C_1 + 6)\exp(-Mnhr_n^2),$$

*where*

$$\begin{aligned}
c_n(M) = \ \ &c_K^2 C_0\sqrt{M}(Mnhr_n^2 + \log n)h^{-1/2}\delta'_n(M)r_n A(h) + \frac{1}{2}c_K C_0 C_1 \sigma^{-2}h^{-1/2}\delta'_n(M)^2 \\
&+ 2c_K^4 C_0^2 C_1 h^{-2}\delta'_n(M)^2(Mnhr_n^2 + \log n)\exp(-(Mnhr_n^2 + \log n)/2).
\end{aligned}$$

*Here $A(h)$ is defined in Section 3.1.*

Based on Theorem 4.1, it is straightforward to extend the rest non-asymptotic theory in Section 3 to the general noise setting. However, we want to point out that the proof of Theorem 4.1 is more involved. The main reason is that now we need to bound higher-order derivatives of the log-likelihood function (which is zero given the quadratic stricture of the log-likelihood function under Gaussian error). Details are deferred to Section S.5 in the supplement for conserving space.

4.2. *Kernel Ridge Regression.* The second extension is to consider the general framework of kernel ridge regression (KRR) [10, 21], where $f_0$ is assumed to belong to a reproducing kernel Hilbert space (RKHS), denoted as $\mathcal{H}$. The corresponding KRR estimator is defined as

$$\widehat{f}_{KRR} = \arg\min_{f \in \mathcal{H}} \left\{ \frac{1}{n} \sum_{i=1}^{n} (Y_i - f(X_i))^2 + \lambda \|f\|_{\mathcal{H}}^2 \right\},$$

where $\|\cdot\|_{\mathcal{H}}$ is the norm associated with $\mathcal{H}$. When $\mathcal{H}$ is chosen as $S^m(\mathbb{I})$, $\widehat{f}_{KRR}$ reduces to the smoothing spline estimate $\widehat{f}_{n,\lambda}$.

We first present a brief review on RKHS theory. Denote $L^2(\mathbb{I})$ as the space of square-integrable functions on $\mathbb{I}$. A subspace of $L^2(\mathbb{I})$ is defined as RKHS if for each $x \in \mathbb{I}$, the evaluation function $f \mapsto f(x)$ is a bounded linear functional. Any RKHS is generated by a positive semidefinite kernel function $K : \mathbb{I} \times \mathbb{I} \to \mathbb{R}$. More precisely, consider the space of all functions of the form

$$g(\cdot) = \sum_{k=1}^{K} \omega_k K(\cdot, v_k), \qquad \text{for some } K \in \mathbb{N}, \ v_k \in \mathbb{I}, \ \omega_k \in \mathbb{R}, \ k = 1, \ldots, K,$$

whose norm is defined by $\|g\|_{\mathcal{H}}^2 := \sum_{k,l=1}^{K} \omega_k \omega_l K(v_k, v_l)$. By taking the closure of this space, it can be shown that we generate an RKHS $\mathcal{H}$ equipped with a norm $\|\cdot\|_{\mathcal{H}}$, and this $\mathcal{H}$ is uniquely associated with the kernel $K$; see [1].

Our finite sample theory can be naturally extended to KRR, by replacing $\widehat{f}_{n,\lambda}$ with $\widehat{f}_{KRR}$ in the test statistic $\widetilde{T}_{n,\lambda}$. After this replacement, it can be shown that the leading term in the separation function $\rho_n^{KRR}(M, L, \lambda)$ for testing the simple hypothesis $H_0 : f = f_0$ becomes

$$(4.3) \qquad \rho_n^{KRR}(M, L, \lambda) \approx \underbrace{\sqrt{\zeta_K \lambda}}_{L_1(\lambda)} + \underbrace{\sqrt{\frac{4 \left( E[K^2(X_1, X_2)] \right)^{1/2}}{n} \sqrt{M}}}_{L_2(\lambda)},$$

where $X_1, X_2 \overset{iid}{\sim} X$. By direct examinations, it can be shown that

$$E(K^2(X_1, X_2)) = \sum_{\nu=1}^{\infty} \frac{1}{(1 + \lambda \rho_\nu)^2},$$

where $\rho_\nu$ are the eigenvalues defined in Section 2. Note that in (4.3), $L_1(\lambda)$ is increasing in $\lambda$ while $L_2(\lambda)$ is decreasing in $\lambda$. To minimize the above leading term, we can select $\lambda = \lambda_{*,KRR}$ by equating $L_1(\lambda)$ and $L_2(\lambda)$, i.e., solving the equation

$$(4.4) \qquad \lambda^{-2} \sum_{\nu=1}^{\infty} \frac{1}{(1 + \lambda \rho_\nu)^2} = \frac{\zeta_K^2 \, n^2}{16 M}.$$

Below is a list of solutions to (4.4) in three concrete situations:

- For polynomial kernel with $\rho_\nu \asymp \nu^{2m}$, the solution of equation (4.4) is

$$\lambda^* \asymp \left( \frac{\sqrt{M}}{\zeta_K \, n} \right)^{\frac{4m}{4m+1}},$$

which recovers our previous result (3.12) by noting the relation that $h = \lambda^{1/(2m)}$.

- For finite rank kernel with $\rho_\nu^{-1} \asymp I(\nu \leq k)$ for a rank $k > 0$, the solution of equation (4.4) is

$$\lambda^* \asymp \frac{\sqrt{k\,M}}{\zeta_K\,n}.$$

- For Gaussian kernel with $\rho_\nu \asymp \exp(\nu^2)$, the solution of equation (4.4) becomes

$$\lambda^* \asymp \frac{\sqrt{M\,\sqrt{\log n}}}{\zeta_K\,n}.$$

**5. Simulation Study.** Simulation results are provided for examining our theory. Consider the following two types of hypotheses:

$$(\text{Simple hypothesis}) \qquad H_0: f = f_0 \text{ versus } H_1: f \neq f_0;$$

$$(\text{Composite hypothesis}) \qquad H_0: f \text{ is linear versus } H_1: f \text{ is not linear},$$

where $f_0 = 5(x^2 - x + \frac{1}{6})$. Set Type I and II errors as $\alpha = \beta = 0.05$.

- For simple hypothesis testing, data were generated as follows

$$(5.1) \quad Y_i = f_c(X_i) + \epsilon_i, \quad \epsilon_i \overset{iid}{\sim} N(0,1), \quad X_i \overset{iid}{\sim} Unif(0,1), \quad \text{and} \quad f_c(x) = \frac{1}{2}c\,x^2 + f_0(x).$$

- For composite hypothesis testing, data were generated as follows

$$(5.2) \quad Y_i = 5X_i + f_c(X_i) + \epsilon_i, \quad \epsilon_i \overset{iid}{\sim} N(0,1), \quad X_i \overset{iid}{\sim} Unif(0,1), \quad \text{and} \quad f_c(x) = c\,(x^2 - x + \frac{1}{6}).$$

Note that the function $f_c$ in model (5.2) lies in the orthogonal complement of the subspace of linear functions in the sense that $E\{f_c(X)\} = E\{Xf_c(X)\} = 0$ for any scaling constant $c$.

We first consider simple hypothesis testing. In some cases, the cutoff value $d_n(\log(15/\alpha), h)$ (see Corollary 3.6) provided in (3.9) can be quite conservative due to the use of some loose concentration inequalities. In practice, we suggest choosing an "exact" cutoff value $d'_n(\log(15/\alpha), h)$ by Monte Carlo simulation. Specifically, by conditioning on $X$, we simulate a number of synthetic datasets $\{\mathbf{Y}^{(k)}\}_{k=1}^N$ from the null model $\mathbf{Y}_i^{(k)} = f_0(X_i) + N(0,1)$, each of which yields a new test statistic $\widetilde{T}_{n,\lambda}^{(k)}$. Set $d'_n(\log(15/\alpha), h)$ as the $(1-\alpha)$-th sample quantile of $\{\widetilde{T}_{n,\lambda}^{(k)}\}_{k=1}^N$.

In simulations, we chose $h$ by directly minimizing $\rho_n(\log(15/\alpha), \log(30/\beta), h)$. Note that we did not replace $d_n$ by $d'_n$ in the above minimization to save computational cost. Such a choice of $h$ is supported by our simulations, and denoted as $h_{FS}$. Then, $d'_n(\log(15/\alpha), h_{FS})$ is used as the cutoff value for the testing procedure. Note that all constants in $\rho_n(M, L, h)$ such as $\rho_K, \zeta_K$ and $c_K$ only depend on $h$ and the eigenvalues of the reproducing kernel operator, which can be well approximately by the empirical eigenvalues of the reproducing kernel matrix.

For simple hypothesis, we compare four testing procedures:

(S1) The proposed $\widetilde{T}_{n,\lambda}$ with the smoothing parameter $h_{FS}$;

(S2) PLRT statistic $PLRT(f_0)$ as described in Remark 3.2 with the same $h_{FS}$;

(S3) The proposed $\widetilde{T}_{n,\lambda}$ with $h$ selected by GCV[2], denoted as $h_{GCV}$;

(S4) PLRT statistic $PLRT(f_0)$ with the same $h_{GCV}$.

The cut-off value for PLRT in S2 and S4 is obtained from Monte Carlo simulation and the null limit distribution given in [20], respectively. Simulation results are reported in Table 1. The rejection proportion (RP) under $c = 0$ reflects the Type I error, while under $c \neq 0$, RP reflects the power.

Overall, all four procedures have comparable type I errors, i.e., $c = 0$, for any sample size. As for power performances, we note that (i) the test using $h_{FS}$ is always more powerful than that using $h_{GCV}$; (ii) $\widehat{T}_{n,\lambda}$ is always more powerful than PLRT given the same choice of $h$. In other words, S1 is always the most powerful one. The observation (i) justifies the finite sample advantage of the non-asymptotic formula in selecting $h$, while (ii) supports the need of removing estimation bias in nonparametric testing; see Remark 3.2. The third observation is that as $c$ increases, $h_{GCV}$ continues decreasing and becomes closer to $h_{FS}$, but never reaches $h_{FS}$. This is consistent with their different asymptotic orders (recall $h_{FS} \asymp n^{-1/(2m+1/2)}$ and $h_{GCV} \asymp n^{-1/(2m+1)}$).

| $n$ | $c$ | $h_{FS}$ | $RP_{FS}$ | $RP_{PLRT}$ | $h_{GCV}$ | $RP_{GCV}$ | $RP'_{GCV}$ |
|---|---|---|---|---|---|---|---|
| 50 | 0 | 0.126 | 0.046 | 0.052 | 0.142(0.008) | 0.052 | 0.054 |
|  | 1 |  | 0.100 | 0.092 | 0.138(0.009) | 0.088 | 0.084 |
|  | 2 |  | 0.396 | 0.340 | 0.134(0.009) | 0.323 | 0.310 |
|  | 3 |  | 0.822 | 0.764 | 0.128(0.009) | 0.752 | 0.748 |
| 100 | 0 | 0.108 | 0.048 | 0.048 | 0.122(0.007) | 0.053 | 0.052 |
|  | 1 |  | 0.264 | 0.170 | 0.117(0.008) | 0.167 | 0.144 |
|  | 2 |  | 0.650 | 0.558 | 0.112(0.007) | 0.493 | 0.474 |
|  | 3 |  | 0.976 | 0.934 | 0.110(0.003) | 0.924 | 0.918 |
| 200 | 0 | 0.092 | 0.050 | 0.048 | 0.104(0.006) | 0.051 | 0.050 |
|  | 1 |  | 0.368 | 0.334 | 0.102(0.007) | 0.325 | 0.290 |
|  | 2 |  | 0.896 | 0.862 | 0.098(0.006) | 0.832 | 0.816 |
|  | 3 |  | 1.00 | 1.00 | 0.094(0.003) | 1.00 | 1.00 |
| 300 | 0 | 0.084 | 0.046 | 0.048 | 0.096(0.006) | 0.048 | 0.048 |
|  | 1 |  | 0.426 | 0.404 | 0.094(0.006) | 0.397 | 0.394 |
|  | 2 |  | 0.968 | 0.946 | 0.092(0.005) | 0.930 | 0.914 |
|  | 3 |  | 1.00 | 1.00 | 0.090(0.005) | 1.00 | 1.00 |
| 400 | 0 | 0.079 | 0.052 | 0.050 | 0.091(0.004) | 0.049 | 0.054 |
|  | 1 |  | 0.668 | 0.640 | 0.087(0.005) | 0.631 | 0.618 |
|  | 2 |  | 1.00 | 1.00 | 0.085(0.004) | 1.00 | 1.00 |
|  | 3 |  | 1.00 | 1.00 | 0.084(0.003) | 1.00 | 1.00 |

TABLE 1

*Simulation results for simple hypothesis testing. $h_{GCV}$ is an average value over 500 replicates (that varies as $c$). $RP_{FS}$, $RP_{PLRT}$, $RP_{GCV}$ and $RP'_{GCV}$ are average rejection proportions by $\widetilde{T}_{n,\lambda}$ with $h_{FS}$, PLRT with $h_{FS}$, $\widehat{T}_{n,\lambda}$ with $h_{GCV}$ and PLRT with $h_{GCV}$ respectively, over 500 replicates.*

For composite hypothesis, we compare two testing procedures:

(C1) The proposed $\widetilde{T}_{n,\lambda}^{com}$ with $h_{FS}^{com}$, selected by numerically minimizing the separation function

---

[2]The $h_{GCV}$ is obtained by using the *ssr* function in the R package *assist*

| $n$ | $c$ | $h_{FS}^{com}$ | $RP_{FS}^{com}$ | $h_{GCV}^{com}$ | $RP_{GCV}^{com}$ |
|-----|-----|------|------|------|------|
| 50  | 0 | 0.247 | 0.060 | 0.310(0.023) | 0.082 |
|     | 1 |       | 0.126 | 0.296(0.020) | 0.084 |
|     | 2 |       | 0.210 | 0.282(0.020) | 0.122 |
|     | 3 |       | 0.348 | 0.270(0.021) | 0.274 |
| 100 | 0 | 0.214 | 0.054 | 0.283(0.020) | 0.064 |
|     | 1 |       | 0.156 | 0.275(0.019) | 0.104 |
|     | 2 |       | 0.234 | 0.265(0.019) | 0.174 |
|     | 3 |       | 0.512 | 0.252(0.016) | 0.380 |
| 200 | 0 | 0.145 | 0.046 | 0.224(0.018) | 0.044 |
|     | 1 |       | 0.220 | 0.208(0.016) | 0.166 |
|     | 2 |       | 0.636 | 0.188(0.015) | 0.536 |
|     | 3 |       | 0.932 | 0.176(0.015) | 0.822 |
| 300 | 0 | 0.102 | 0.054 | 0.184(0.017) | 0.046 |
|     | 1 |       | 0.254 | 0.175(0.013) | 0.190 |
|     | 2 |       | 0.714 | 0.165(0.015) | 0.656 |
|     | 3 |       | 0.976 | 0.153(0.014) | 0.882 |
| 400 | 0 | 0.096 | 0.054 | 0.164(0.013) | 0.050 |
|     | 1 |       | 0.290 | 0.156(0.014) | 0.256 |
|     | 2 |       | 0.862 | 0.146(0.013) | 0.788 |
|     | 3 |       | 1.00  | 0.138(0.012) | 0.946 |

TABLE 2

*Simulation results for composite hypothesis testing. $h_{GCV}^{com}$ is an average value over 500 replicates (that varies as $c$). $RP_{FS}^{com}$ and $RP_{GCV}^{com}$ are average rejection proportions by $\overline{T}_{n,\lambda}^{com}$ with $h_{FS}^{com}$ and PLRT with $h_{GCV}^{com}$, respectively, over 500 replicates.*

in Theorem 3.7;

(C2) PLRT statistics $PLRT(f_0)$ with $h$ selected by GCV, denoted as $h_{GCV}$;.

Simulation results are reported in Table 2. We observe similar phenomena as in the simple testing case. In particular, C1 is uniformly more powerful than C2 due to the use of $h_{FS}^{com}$.

**6. Proofs.** In this section, we provide proofs of the main theorems in the paper. More technical details are provided in the supplementary material.

6.1. *Preliminaries for penalized likelihood estimation.* Before formal proofs, let us introduce some preliminaries. Considering model (2.1) parametrized by the unknown regression function $f \in H^m(1)$. The Fréchet derivative of the penalized loss function $\ell_{n,\lambda}$ at $g \in S^m(\mathbb{I})$ can be identified as

$$D\ell_{n,\lambda}(g)g_1 = \frac{1}{n}\sum_{i=1}^{n}\ell^{(1)}(Y_i - g(X_i))g_1(X_i) - \langle \mathcal{P}_\lambda g, g_1 \rangle \equiv \langle S_{n,\lambda}(g), g_1 \rangle,$$

when operated on arbitrary function $g_1$ in $S^m(\mathbb{I})$. Let $S_\lambda(g) = E_f\{S_{n,\lambda}(g)\}$ be expectation of $S_{n,\lambda}(g)$, for any $g \in S^m(\mathbb{I})$. We denote the second- and third-order Fréchet derivatives of $S_\lambda$ by $DS_\lambda$ and $D^2 S_\lambda$. By the optimality of the smoothing spline estimator $\widehat{f}_{n,\lambda}$, we have $S_{n,\lambda}(\widehat{f}_{n,\lambda}) = 0$. Therefore, $S_{n,\lambda}(f)$ can be expressed as

$$(6.1) \qquad S_{n,\lambda}(f) = \frac{1}{n} \sum_{i=1}^{n} \epsilon_i K_{X_i} - \mathcal{P}_\lambda f.$$

The Fréchet derivatives of $S_{n,\lambda}$ and $DS_{n,\lambda}$ are denoted $DS_{n,\lambda}(g)g_1 g_2$ and $D^2 S_{n,\lambda}(g)g_1 g_2 g_3$. These derivatives can be explicitly written as

$$D^2 \ell_{n,\lambda}(g)g_1 g_2 \equiv DS_{n,\lambda}(g)g_1 g_2 = n^{-1} \sum_{i=1}^{n} \ell^{(2)}(Y_i - g(X_i))g_1(X_i)g_2(X_i) - \langle \mathcal{P}_\lambda g_1, g_2 \rangle,$$

$$D^3 \ell_{n,\lambda}(g)g_1 g_2 g_3 \equiv D^2 S_{n,\lambda}(g)g_1 g_2 g_3 = n^{-1} \sum_{i=1}^{n} \ell^{(3)}(Y_i - g(X_i))g_1(X_i)g_2(X_i)g_3(X_i).$$

We have the following concentration inequality for sum of iid elements $\{K_{X_i} : i = 1, \ldots, n\}$, which plays a critical role in proving our large deviation bounds for the smoothing spline estimator. A proof is deferred to Section 6.5. Recall that $\mathcal{G} := \{g \in S^m(\mathbb{I}) : \|g\|_{\sup} \leq 1, J(g,g) \leq c_K^{-2} h^{-2m+1}\}$.

LEMMA 6.1.    *Suppose that $\psi_{n,f}(z; g)$ is a measurable function defined upon $z = (y,x) \in \mathcal{Y} \times \mathbb{I}$ and $g \in \mathcal{G}$ satisfying $\psi_{n,f}(z; 0) = 0$ and the following Lipschitz continuity condition: for any $f \in H^m(1)$, $i = 1, \ldots, n$ and any $g_1, g_2 \in \mathcal{G}$,*

$$(6.2) \qquad |\psi_{n,f}(Z_i; g_1) - \psi_{n,f}(Z_i; g_2)| \leq c_K^{-1} h^{1/2} \|g_1 - g_2\|_{\sup}.$$

*Then for any constant $t \geq 0$ and $n \geq 1$,*

$$\sup_{f \in H^m(1)} P_f^n \left( \sup_{g \in \mathcal{G}} \|Z_{n,f}(g)\| > t \right) \leq 2 \exp\left( -\frac{t^2}{A(h)^2} \right),$$

*where recall $A(h) = A(h, 2)$ and*

$$Z_{n,f}(g) = \frac{1}{\sqrt{n}} \sum_{i=1}^{n} [\psi_{n,f}(Z_i; g)K_{X_i} - E_f\{\psi_{n,f}(Z_i; g)K_{X_i}\}].$$

6.2. *Proofs of large deviation bounds for smoothing spline estimates.*    Given the development in the previous subsection, we are now ready to prove the two large deviation bounds in Theorem 3.1 and Theorem 3.4.

*Proof of Theorem 3.1.*    For any $g, g_1 \in S^m(\mathbb{I})$, we have the following sequence of identities,

$$\ell_{n,\lambda}(g) = -\frac{1}{2n} \sum_{i=1}^{n} (Y_i - g(X_i))^2 - \frac{\lambda}{2} J(g),$$

$$S_{n,\lambda}(g) = \frac{1}{n} \sum_{i=1}^{n} (Y_i - g(X_i))K_{X_i} - \mathcal{P}_\lambda g,$$

$$DS_{n,\lambda}(g)g_1 = -\frac{1}{n} \sum_{i=1}^{n} g_1(X_i)K_{X_i} - \mathcal{P}_\lambda g_1, \quad \text{and} \quad D^2 S_{n,\lambda}(g) = 0.$$

For any $f \in S^m(\mathbb{I})$, define $f_\lambda = f + S_\lambda(f)$. It follows from [20, Proposition 2.3] that $DS_\lambda(f) = -id$, where $id$ denotes the identity operator. Then

$$S_\lambda(f_\lambda) - S_\lambda(f) = DS_\lambda(f)(f_\lambda - f) = -(f_\lambda - f) = -S_\lambda(f),$$

hence, $S_\lambda(f_\lambda) = 0$. Meanwhile,

$$\|f_\lambda - f\| = \|S_\lambda(f)\| = \|\mathcal{P}_\lambda f\| \le \sqrt{\lambda J(f)} \le h^m.$$

By using the property $S_\lambda(f_\lambda) = 0$, we have the following bound

$$
\begin{aligned}
\|S_{n,\lambda}(f_\lambda)\| &= \|S_{n,\lambda}(f_\lambda) - S_\lambda(f_\lambda)\| \\
&\le \|\frac{1}{n}\sum_{i=1}^n [(f - f_\lambda)(X_i)K_{X_i} - E\{(f - f_\lambda)(X)K_X\}]\| + \|\frac{1}{n}\sum_{i=1}^n \epsilon_i K_{X_i}\|.
\end{aligned}
$$

Now we bound the two terms on the right hand side, respectively.

To bound the first term, we apply Lemma 6.1, with $\psi_{n,f}(Z; g) = g/(c_K h^{-1/2}) = c_K^{-1} h^{1/2} g(X)$, to obtain

$$(6.3) \quad P\left(\sup_{\bar{g} \in \mathcal{G}} \|\frac{1}{n}\sum_{i=1}^n [\bar{g}(X_i)K_{X_i} - E\{\bar{g}(X)K_X\}]\| > c_K(nh)^{-1/2}t\right) \le \exp(-t^2/A(h)^2), \ t > 0.$$

If we let the event

$$\mathcal{E}_n := \left\{\sup_{\bar{g} \in \mathcal{G}} \|\frac{1}{n}\sum_{i=1}^n [\bar{g}(X_i)K_{X_i} - E\{\bar{g}(X)K_X\}]\| \le c_K M^{1/2} r A(h)\right\},$$

then (6.3) implies $P(\mathcal{E}_n^c) \le \exp(-Mnhr^2)$. Since we have $\bar{f} = (f - f_\lambda)/(c_K h^{-1/2}\|f - f_\lambda\|) \in \mathcal{G}$, we obtain that on $\mathcal{E}_n$,

$$\|\frac{1}{n}\sum_{i=1}^n [\bar{f}(X_i)K_{X_i} - E\{\bar{f}(X)K_X\}]\| \le c_K M^{1/2} r A(h),$$

leading to the following bound for the first term

$$\|\frac{1}{n}\sum_{i=1}^n [(f - f_\lambda)(X_i)K_{X_i} - E\{(f - f_\lambda)(X)K_X\}]\| \le c_K^2 M^{1/2} h^{m-1/2} r A(h).$$

To bound the second term, let $\Sigma = [K(X_i, X_j)]_{1 \le i,j \le n}$ and $\boldsymbol{\epsilon} = (\epsilon_1, \ldots, \epsilon_n)^T$. By the Hanson-Wright inequality (see, for example, [11]), we have

$$P(\boldsymbol{\epsilon}^T \Sigma \boldsymbol{\epsilon} > tr(\Sigma) + 2\sqrt{tr(\Sigma^2)Mnhr^2} + 2\|\Sigma\|_F Mnhr^2) \le \exp(-Mnhr^2),$$

where $\|\cdot\|_F$ is the Frobenius norm. Since

$$
\begin{aligned}
tr(\Sigma) &= \sum_{i=1}^n K(X_i, X_i) \le c_K^2 nh^{-1}, \\
tr(\Sigma^2) &= \sum_{i,j=1}^n K(X_i, X_j)^2 \le c_K^4 n^2 h^{-2}, \\
\|\Sigma\|_F &\le \sqrt{tr(\Sigma^2)} \le c_K^2 nh^{-1},
\end{aligned}
$$

we get that

$$P\left(\|\frac{1}{n}\sum_{i=1}^{n}\epsilon_i K_{X_i}\| > c_K(\sqrt{2M}r + (nh)^{-1/2})\right) \leq \exp(-Mnhr^2).$$

Now if we let event

$$\mathcal{E}_n' = \{\|\frac{1}{n}\sum_{i=1}^{n}\epsilon_i K_{X_i}\| \leq c_K(\sqrt{2M}r + (nh)^{-1/2})\},$$

then we have $P(\mathcal{E}_n') \geq 1 - \exp(-Mnhr^2)$, and on $\mathcal{E}_n'$, the second term can be bounded as

$$\|\frac{1}{n}\sum_{i=1}^{n}\epsilon_i K_{X_i}\| \leq c_K(\sqrt{2M}r + (nh)^{-1/2}).$$

Putting pieces together, we obtain that with probability at least $P(\mathcal{E}_n \cap \mathcal{E}_n') \geq 1 - 2\exp(-Mnhr^2)$,

$$(6.4) \qquad \|S_{n,\lambda}(f_\lambda)\| \ \leq \ c_K^2 M^{1/2} h^{m-1/2} r A(h) + c_K(\sqrt{2M}r + (nh)^{-1/2}) := r'/2.$$

Now let us consider the following class of operators, indexed by functions $f$ in $S^m(\mathbb{I})$, as

$$T_f(g) = g + S_{n,\lambda}(f_\lambda + g), \quad \text{for all } g \in S^m(\mathbb{I}).$$

By adding and subtracting the same term, we can express $T_f$ as

$$\begin{aligned}
T_f(g) \ &= \ -DS_\lambda(f_\lambda)^{-1}[DS_{n,\lambda}(f_\lambda)g - DS_\lambda(f_\lambda)g] \\
&\quad -DS_\lambda(f_\lambda)^{-1}[S_{n,\lambda}(f_\lambda + g) - S_{n,\lambda}(f_\lambda) - DS_{n,\lambda}(f_\lambda)g] \\
&\quad -DS_\lambda(f_\lambda)^{-1}S_{n,\lambda}(f_\lambda) \\
&= \ DS_{n,\lambda}(f_\lambda)g - DS_\lambda(f_\lambda)g + S_{n,\lambda}(f_\lambda) \\
&= \ -\frac{1}{n}\sum_{i=1}^{n}[g(X_i)K_{X_i} - E\{g(X)K_X\}] + S_{n,\lambda}(f_\lambda).
\end{aligned}$$

By (6.4), we obtain that with probability at least $1 - 2\exp(-Mnhr^2)$, for any $g \in S^m(\mathbb{I})$,

$$(6.5) \qquad \|T_f(g)\| \leq \|\frac{1}{n}\sum_{i=1}^{n}[g(X_i)K_{X_i} - E\{g(X)K_X\}]\| + r'/2$$

$$(6.6) \qquad \leq c_K^2 M^{1/2} rh^{-1/2} A(h)\|g\| + r'/2 \leq r',$$

where the last inequality follows by the condition $c_K^2 M^{1/2} rh^{-1/2} A(h) \leq 1/2$. Since for any $g_1, g_2 \in S^m(\mathbb{I})$,

$$\begin{aligned}
&\|T_f(g_1) - T_f(g_2)\| \\
= \ &\|\frac{1}{n}\sum_{i=1}^{n}[(g_1 - g_2)(X_i)K_{X_i} - E\{(g_1 - g_2)(X)K_X\}]\| \\
\leq \ &c_K^2 M^{1/2} rh^{-1/2} A(h)\|g_1 - g_2\| \leq (1/2)\|g_1 - g_2\|,
\end{aligned}$$

$T_f$ is a contraction mapping on $\mathbb{B}(r')$. Therefore, by the contraction mapping theorem ([19]), there exists $g \in \mathbb{B}(r')$ such that $S_{n,\lambda}(f_\lambda + g) = 0$. Since the smoothing spline estimate $\widehat{f}_{n,\lambda}$ also satisfies $S_{n\lambda}(\widehat{f}_{n,\lambda}) = 0$, and the solution is unique, we must have $\widehat{f}_{n,\lambda} = f_\lambda + g$. Hence, inequality (6.5) yields $\|\widehat{f}_{n,\lambda} - f\| = \|g\| \le h^m + r'$. Consequently, we obtain

$$\sup_{f \in H^m(1)} P_f\left(\|\widehat{f}_{n,\lambda} - f\| \ge \delta_n(M,r)\right) \le 2\exp(-Mnhr^2),$$

where $\delta_n(M,r) = 2h^m + c_K(\sqrt{2M}r + (nh)^{-1/2})$. This expression of $\delta_n$ follows from $c_K^2 M^{1/2} rh^{-1/2} A(h) \le 1/2$, leading to $r' = 2c_K^2 M^{1/2} h^{m-1/2} rA(h) + 2c_K(\sqrt{2M}r + (nh)^{-1/2}) \le h^m + 2c_K(\sqrt{2M}r + (nh)^{-1/2})$.

*Proof of Theorem 3.4.* Given the development in the proof of Theorem 3.1, the proof for this theorem is easy. Let us define $g_n = \widehat{f}_{n,\lambda} - f$. Note that on $\mathcal{E}_n \cap \mathcal{E}'_n$

$$
\begin{aligned}
&\|S_{n,\lambda}(f + g_n) - S_{n,\lambda}(f) - (S_\lambda(f + g_n) - S_\lambda(f))\| \\
=\ &\|\frac{1}{n}\sum_{i=1}^n [g_n(X_i)K_{X_i} - E\{g_n(X)K_X\}]\| \\
\le\ &c_K^2 M^{1/2} rh^{-1/2} A(h)\|g_n\| \\
\le\ &c_K^2 M^{1/2} rh^{-1/2} A(h)\delta_n(M,r).
\end{aligned}
$$

Moreover, we have the identity

$$
\begin{aligned}
&S_{n,\lambda}(f + g_n) - S_{n,\lambda}(f) - (S_\lambda(f + g_n) - S_\lambda(f)) \\
=\ &0 - S_{n,\lambda}(f) - DS_\lambda(f)g_n = g_n - S_{n,\lambda}(f).
\end{aligned}
$$

Therefore, we have

$$\sup_{f \in H^m(1)} P_f\left(\|\widehat{f}_{n,\lambda} - f - S_{n,\lambda}(f)\| \ge c_K^2 M^{1/2} rh^{-1/2} A(h)\delta_n(M,r)\right) \le 2\exp(-Mnhr^2).$$

### 6.3. *Proof of Theorem 3.5.*
We will make use of the following two lemmas that control the type I and type II errors of the test. Their proofs are provided in the supplement.

LEMMA 6.2 (Type I error). *For any $M \ge 0$ satisfying the conditions of Theorem 3.4, it holds that*

$$(6.7) \qquad P_{f_0}\left(|\widetilde{T}_{n,\lambda}| \ge \frac{4\rho_K}{n\sqrt{h}}\sqrt{M} + R_{1,n}(M)\right) \le 15\exp(-M),$$

*where $\rho_K^2 = hE[K^2(X_1, X_2)]$ with $X_1, X_2 \overset{iid}{\sim} X$, and the explicit form of the remainder term $R_{1,n}(M)$ is provided in Section S.2.*

LEMMA 6.3 (Type II error). *For any $L > 0$ satisfying the conditions of Theorem 3.4, we have*

$$
\sup_{f-f_0 \in H^m(1)} P_f\left(|\widetilde{T}_{n,\lambda} - \|(I-\mathcal{P}_\lambda)(f-f_0)\|^2| \geq \left(\frac{2\sqrt{2}\,\|(I-\mathcal{P}_\lambda)(f-f_0)\|}{\sqrt{n}} + \frac{4\rho_K}{n\sqrt{h}}\right)\sqrt{L}\right.
$$
$$
\left. + R_{2,n}(L)\right) \leq 30\exp(-L),
$$
(6.8)

*where the explicit form of the remainder term $R_{2,n}(L)$ is provided in Section S.2.*

Now, we will prove the theorem. By applying Lemma 6.2, if we set $M = \log(15\,\alpha^{-1})$, then

$$
P_{f_0}(|\widetilde{T}_{n,\lambda}| > d_n(M,h)) \leq \alpha.
$$

By using Lemma 6.3, we obtain

$$
\sup_{\|f-f_0\| \geq \rho_n(M,L,h)} P_f(|\widetilde{T}_{n,\lambda}| < d_n(M,h))
$$
$$
= \sup_{\|f-f_0\| \geq \rho_n(M,L,h)} P_f\left(\|(I-\mathcal{P}_\lambda)(f-f_0)\|^2 - |\widetilde{T}_{n,\lambda}| \geq -d_n(M,h) + \left(\|(I-\mathcal{P}_\lambda)(f-f_0)\| - \sqrt{\frac{2L}{n}}\right)^2\right.
$$
$$
\left. + \frac{2\sqrt{2}\,\|(I-\mathcal{P}_\lambda)(f-f_0)\|}{\sqrt{n}}\sqrt{L} - \frac{2L}{n}\right)
$$
$$
\leq \sup_{\|f-f_0\| \geq \rho_n(M,L,h)} P_f\left(|\widetilde{T}_{n,\lambda}| - \|(I-\mathcal{P}_\lambda)(f-f_0)\|^2| \geq \frac{2\sqrt{2}\,\|(I-\mathcal{P}_\lambda)(f-f_0)\|}{\sqrt{n}}\sqrt{L}\right.
$$
$$
\left. + \left(\rho_n(M,L,h) - \sqrt{\zeta_K\lambda} - \sqrt{\frac{2L}{n}}\right)^2 - d_n(M,h) - \frac{2L}{n}\right)
$$
$$
\leq 30\exp(-L),
$$

where $\zeta_K = \sup_{g \in H^m(1)} \lambda^{-1}\|\mathcal{P}_\lambda g\|^2$ and $\rho_n(M,L,h)$ is the solution to the equation:

$$
\left(\rho_n(M,L,h) - \sqrt{\zeta_K\lambda} - \sqrt{\frac{2L}{n}}\right)^2 - d_n(M,h) - \frac{2L}{n} = R_{2,n}(L),
$$
$$
\text{i.e.} \quad \rho_n(M,L,h) = \sqrt{\zeta_K\lambda} + \sqrt{\frac{2L}{n}} + \sqrt{d_n(M) + \frac{2L}{n} + R_{2,n}(L)}.
$$

6.4. *Proof of Theorem 3.7.* Similar to the simple hypothesis testing case, we make use of the following two lemmas that control the type I and type II errors of the composite test. Their proofs are provided in the supplement.

LEMMA 6.4 (Type I error). *If $f_0$ is a linear function and $M \in (0, n/4]$ satisfies the conditions of Theorem 3.4, then it holds that*

$$
P_{f_0}\left(|\widetilde{T}_{n,\lambda}^{com}| \geq \frac{4\rho_K}{n\sqrt{h}}\sqrt{M} + R_{1,n}^c(M)\right) \leq 24\exp(-M),
$$

*where the form of $R_{1,n}^c(M)$ is provided in Section S.2.*

LEMMA 6.5 (Type II error). *For any $L > 0$ satisfying the conditions of Theorem 3.4, it holds that*

$$\sup_{f \in H^m(1)} P_f^n \Big( \big| \widetilde{T}_{n,\lambda}^{com} - \|f^\perp - \mathcal{P}_\lambda f^\perp\|^2 \big| \geq \Big( \frac{2\sqrt{2}\,\|f^\perp - \mathcal{P}_\lambda f^\perp\|}{\sqrt{n}} + \frac{4\rho_K}{n\sqrt{h}} \Big) \sqrt{L}$$
$$+ R_{2,n}^c(L) \Big) \leq 60 \exp(-L),$$

*where the form of $R_{2,n}^c(L)$ is provided in Section S.2.*

Now we prove the theorem. By applying Lemma 6.4, if we set $M = \log(24\,\alpha^{-1})$, then

$$\sup_{\substack{f_0 \in H^m(1) \\ f_0 \text{ is linear}}} P_{f_0}(|\widetilde{T}_{n,\lambda}^{com}| > d_n^{com}(M, h)) \leq \alpha.$$

By using Lemma 6.5, we obtain

$$\sup_{\substack{f \in H^m(1) \\ \|f^\perp\| \geq \rho_n^{com}(M, L, h)}} P_f(|\widetilde{T}_{n,\lambda}^{com}| < d_n^{com}(M, h))$$

$$= \sup_{\substack{f \in H^m(1) \\ \|f^\perp\| \geq \rho_n^{com}(M, L, h)}} P_f\Big( \|f^\perp - \mathcal{P}_\lambda f^\perp\|^2 - |\widetilde{T}_{n,\lambda}^{com}| \geq -d_n(M, h) + \Big( \|f^\perp - \mathcal{P}_\lambda f^\perp\| - \sqrt{\frac{2L}{n}} \Big)^2$$
$$+ \frac{2\sqrt{2}\,\|f^\perp - \mathcal{P}_\lambda f^\perp\|}{\sqrt{n}} \sqrt{L} - \frac{2L}{n} \Big)$$

$$\leq \sup_{\substack{f \in H^m(1) \\ \|f^\perp\| \geq \rho_n^{com}(M, L, h)}} P_f\Big( \big| |\widetilde{T}_{n,\lambda}^{com}| - \|f^\perp - \mathcal{P}_\lambda f^\perp\| \big| \geq \frac{2\sqrt{2}\,\|f^\perp - \mathcal{P}_\lambda f^\perp\|}{\sqrt{n}} \sqrt{L}$$
$$+ \Big( \rho_n^{com}(M, L, h) - \sqrt{\zeta_K\lambda} - \sqrt{\frac{2L}{n}} \Big)^2 - d_n^{com}(M, h) - \frac{2L}{n} \Big)$$

$$\leq (2e + 26)\exp(-L),$$

where $\zeta_K = \sup_{g \in H^m(1)} \lambda^{-1} \|\mathcal{P}_\lambda g\|^2$ and $\rho_n^{com}(M, L, h)$ is the solution to the equation:

$$\Big( \rho_n^{com}(M, L, h) - \sqrt{\zeta_K\lambda} - \sqrt{\frac{2L}{n}} \Big)^2 - d_n^{com}(M, h) - \frac{2L}{n} = R_{2,n}^c(L),$$

i.e. $\rho_n^{com}(M, L, h) = \sqrt{\zeta_K\lambda} + \sqrt{\frac{2L}{n}} + \sqrt{d_n^{com}(M, h) + \frac{2L}{n} + R_{2,n}^c(L)}.$

6.5. *Proof of Lemma 6.1.* For any $f \in H^m(1)$ and $n \geq 1$, and any functions $g_1, g_2 \in \mathcal{G}$, we have the following bound for each additive component in the sum,

$$\|(\psi_{n,f}(Z_i; g_1) - \psi_{n,f}(Z_i; g_2))K_{X_i}\| \leq c_K^{-1} h^{1/2} \|g_1 - g_2\|_{\sup} c_K h^{-1/2} = \|g_1 - g_2\|_{\sup}.$$

For fixed $g_1, g_2$, we apply the bounded difference inequality (see, for example, Theorem 3.5 of [18]) to obtain the following concentration inequality for the sum,

$$(6.9) \qquad P_f^n\left( \|Z_{n,f}(g_1) - Z_{n,f}(g_2)\| \geq t \right) \leq 2\exp\left( -\frac{t^2}{8\|g_1 - g_2\|_{\sup}^2} \right), \quad \text{for any } t > 0.$$

Due to the equivalence between the sub-Gaussian tail of a random variable and its Orlicz norm (Lemma 8.1 in [13]), we obtain

$$\Big\| \|Z_{n,f}(g_1) - Z_{n,f}(g_2)\| \Big\|_{\psi_2} \leq \sqrt{24} \, \|g_1 - g_2\|_{\sup},$$

where $\| \cdot \|_{\psi_2}$ denotes the Orlicz norm associated with $\psi_2(s) \equiv \exp(s^2) - 1$. Let $\tau = \sqrt{\log 1.5} \approx 0.6368$ and $\phi(M) = \psi_2(\tau x)$. Then it is easy to check that $\phi(1) \leq 1/2$, and for any $x, y \geq 1$, $\phi(M) \, \phi(y) \leq \phi(xy)$. By applying Lemma 8.2 in [13], we have the following relationship between the Orlicz norm of the max of $l$ ($l \in \mathbb{N}$) random variables and the max of their individual Orlicz norms,

$$(6.10) \qquad \| \max_{1 \leq i \leq l} \xi_i \|_{\psi_2} \leq \frac{2}{\tau} \psi_2^{-1}(l) \max_{1 \leq i \leq l} \|\xi_i\|_{\psi_2},$$

for any random variables $\xi_1, \ldots, \xi_l$.

Next we apply a "chaining" argument to prove the desired concentration inequality based on (6.9). Let $T_0 \subset T_1 \subset T_2 \subset \cdots \subset T_\infty \equiv \mathcal{G}$ be a sequence of finite nested sets satisfying the following properties:

(i). for any $T_j$ and any $s, t \in T_j$, $\|s - t\|_{\sup} \geq \varepsilon 2^{-j}$; each $T_j$ is "maximal" in the sense that if one adds any point in $T_j$, then the inequality will fail;

(ii). the cardinality of $T_j$ is upper bounded by

$$\log |T_j| \leq \log N(\varepsilon 2^{-j}, \mathcal{G}, \| \cdot \|_{\sup}) \leq c_0 c_K^{-1/m} h^{-(2m-1)/(2m)} (\varepsilon 2^{-j})^{-1/m},$$

where $c_0 > 0$ is absolute constant; Here we used the fact that the covering entropy of an $R$-ball in $S^m(\mathbb{I})$ relative to the sup-norm is of order $R^{-1/m}$.

(iii). each element $t_{j+1} \in T_{j+1}$ is uniquely linked to an element $t_j \in T_j$ which satisfies $\|t_j - t_{j+1}\|_{\sup} \leq \varepsilon 2^{-j}$.

Based on this sequence $\{T_k : k \geq 0\}$, for arbitrary $s_{k+1}, t_{k+1} \in T_{k+1}$ with $\|s_{k+1} - t_{k+1}\|_{\sup} \leq \varepsilon$ we can choose two chains (both of length $k + 2$) $\{t_j : 0 \leq j \leq k+1\}$ and $\{s_j : 0 \leq j \leq k+1\}$ with $t_j, s_j \in T_j$ for $0 \leq j \leq k+1$, such that the end points $s_0$ and $t_0$ satisfy

$$\|s_0 - t_0\|_{\sup} \leq \sum_{j=0}^{k} [\|s_j - s_{j+1}\|_{\sup} + \|t_j - t_{j+1}\|_{\sup}] + \|s_{k+1} - t_{k+1}\|_{\sup} \leq 2 \sum_{j=0}^{k} \varepsilon 2^{-j} + \varepsilon \leq 5\varepsilon,$$

implying $\Big\| \|Z_{n,f}(s_0) - Z_{n,f}(t_0)\| \Big\|_{\psi_2} \leq 5\sqrt{24} \, \varepsilon$. Recall that function $\Psi(r) = \int_0^r \sqrt{\log(1 + \exp(x^{-1/m}))} \, dx$.

Now, it follows from (6.10) that

$$
\left\| \max_{s_{k+1}, t_{k+1} \in T_{k+1}} \| Z_{n,f}(s_{k+1}) - Z_{n,f}(t_{k+1}) - (Z_{n,f}(s_0) - Z_{n,f}(t_0)) \| \right\|_{\psi_2}
$$

$$
\leq 2 \sum_{j=0}^{k} \left\| \max_{\substack{u \in T_{j+1}, v \in T_j \\ u, v \text{ link each other}}} \| Z_{n,f}(u) - Z_{n,f}(v) \| \right\|_{\psi_2}
$$

$$
\leq \frac{4}{\tau} \sum_{j=0}^{k} \psi_2^{-1}(N(2^{-j-1}\varepsilon, \mathcal{G}, \| \cdot \|_{\sup})) \max_{\substack{u \in T_{j+1}, v \in T_j \\ u, v \text{ link each other}}} \left\| \| Z_{n,f}(u) - Z_{n,f}(v) \| \right\|_{\psi_2}
$$

$$
\leq \frac{4\sqrt{24}}{\tau} \sum_{j=0}^{k} \sqrt{\log\left(1 + N(\varepsilon 2^{-j-1}, \mathcal{G}, \| \cdot \|_{\sup})\right)} \, \varepsilon \, 2^{-j}
$$

$$
\leq \frac{8\sqrt{24}}{\tau} \sum_{j=1}^{k+1} \sqrt{\log\left(1 + \exp\left(c_0 c_K^{-1/m} h^{-(2m-1)/(2m)}(\varepsilon 2^{-j})^{-1/m}\right)\right)} \, \varepsilon \, 2^{-j}
$$

$$
\leq \frac{32\sqrt{6}}{\tau} \int_0^{\varepsilon/2} \sqrt{\log\left(1 + \exp\left(c_0 c_K^{-1/m} h^{-(2m-1)/(2m)} M^{-1/m}\right)\right)} \, dx
$$

$$
= \frac{32\sqrt{6}}{\tau} c_K^{-1} c_0^m h^{-(2m-1)/2} \Psi\left(\frac{1}{2} c_K c_0^{-m} h^{(2m-1)/2} \varepsilon\right).
$$

On the other hand, we have

$$
\left\| \max_{\substack{u,v \in T_0 \\ \|u-v\|_{\sup} \leq 5\varepsilon}} \| Z_{n,f}(u) - Z_{n,f}(v) \| \right\|_{\psi_2}
$$

$$
\leq \frac{2}{\tau} \psi_2(|T_0|^2) \max_{\substack{u,v \in T_0 \\ \|u-v\|_{\sup} \leq 5\varepsilon}} \| \| Z_{n,f}(u) - Z_{n,f}(v) \| \|_{\psi_2} \leq \frac{2}{\tau} \psi_2^{-1}(N(\varepsilon, \mathcal{G}, \| \cdot \|_{\sup})^2)(5\sqrt{24}\varepsilon).
$$

Combining these two bounds together, we obtain

$$
\left\| \max_{\substack{s,t \in T_{k+1} \\ \|s-t\|_{\sup} \leq \varepsilon}} \| Z_{n,f}(s) - Z_{n,f}(t) \| \right\|_{\psi_2}
$$

$$
\leq \left\| \max_{s_{k+1}, t_{k+1} \in T_{k+1}} \| Z_{n,f}(s_{k+1}) - Z_{n,f}(t_{k+1}) - (Z_{n,f}(s_0) - Z_{n,f}(t_0)) \| \right\|_{\psi_2}
$$

$$
+ \left\| \max_{\substack{u,v \in T_0 \\ \|u-v\|_{\sup} \leq 5\varepsilon}} \| Z_{n,f}(u) - Z_{n,f}(v) \| \right\|_{\psi_2}
$$

$$
\leq \frac{32\sqrt{6}}{\tau} c_K^{-1} c_0^m h^{-(2m-1)/2} \Psi\left(\frac{1}{2} c_K c_0^{-m} h^{(2m-1)/2} \varepsilon\right) + \frac{2}{\tau} \psi_2^{-1}(N(\varepsilon, \mathcal{G}, \| \cdot \|_{\sup})^2)(5\sqrt{24}\varepsilon)
$$

$$
\leq \frac{32\sqrt{6}}{\tau} c_K^{-1} c_0^m h^{-(2m-1)/2} \Psi\left(\frac{1}{2} c_K c_0^{-m} h^{(2m-1)/2} \varepsilon\right)
$$

$$
+ \frac{10\sqrt{24}\varepsilon}{\tau} \sqrt{\log\left(1 + \exp\left(2 c_0 (c_K h^{(2m-1)/2} \varepsilon)^{-1/m}\right)\right)}
$$

$$
= A(h, \varepsilon),
$$

where in the last step we used the definition of $A(h, \varepsilon)$.

Now consider any two functions $g_1, g_2 \in \mathcal{G}$ with $\|g_1 - g_2\|_{\sup} \leq \varepsilon/2$. For any $k \geq 2$, since $T_k$ is "maximal" due to our construction, there must exist $(s_k, t_k) \in T_k^2$ such that $\max\{\|g_1 - s_k\|_{\sup}, \|g_2 - t_k\|_{\sup}\} \leq \varepsilon 2^{-k}$, which also implies $\|s_k - t_k\|_{\sup} \leq \varepsilon$. Therefore, we can decompose the difference between $Z_{n,f}(g_1)$ and $Z_{n,f}(g_2)$ by

$$
\begin{aligned}
&\|Z_{n,f}(g_1) - Z_{n,f}(g_2)\| \\
\leq \ & \|Z_{n,f}(g_1) - Z_{n,f}(s_k)\| + \|Z_{n,f}(g_2) - Z_{n,f}(t_k)\| + \|Z_{n,f}(s_k) - Z_{n,f}(t_k)\| \\
\leq \ & 4\sqrt{n}\,\varepsilon\, 2^{-k} + \max_{\substack{u,v \in T_k \\ \|u-v\|_{\sup} \leq \varepsilon}} \|Z_{n,f}(u) - Z_{n,f}(v)\|.
\end{aligned}
$$

Now we can obtain

$$
\begin{aligned}
& \left\| \sup_{\substack{g_1, g_2 \in \mathcal{G} \\ \|g_1-g_2\|_{\sup} \leq \varepsilon/2}} \|Z_{n,f}(g_1) - Z_{n,f}(g_2)\| \right\|_{\psi_2} \\
\leq \ & 4\sqrt{n}\varepsilon 2^{-k}/\sqrt{\log 2} + \left\| \max_{\substack{u,v \in T_k \\ \|u-v\|_{\sup} \leq \varepsilon}} \|Z_{n,f}(u) - Z_{n,f}(v)\| \right\|_{\psi_2} \\
\leq \ & 4\sqrt{n}\varepsilon 2^{-k}/\sqrt{\log 2} + A(h, \varepsilon) \to A(h, \varepsilon), \quad \text{by letting } k \to \infty.
\end{aligned}
$$

Taking $\varepsilon = 2$ in the above inequality, we obtain

$$
\left\| \sup_{\substack{g_1, g_2 \in \mathcal{G} \\ \|g_1-g_2\|_{\sup} \leq 1}} \|Z_{n,f}(g_1) - Z_{n,f}(g_2)\| \right\|_{\psi_2} \leq A(h, 2) = A(h).
$$

Consequently, by choosing $g_2 = 0$ and using the equivalence between sub-Gaussian tails and the Orlicz norm (Lemma 8.1 in [13]), we obtain

$$
P_f^n\left( \sup_{g \in \mathcal{G}} \|Z_{n,f}(g)\| \geq t \right) \leq 2\exp\left(-\frac{t^2}{A(h)^2}\right).
$$

This completes the proof.

*Supplementary document to*

Non-asymptotic Theory for Nonparametric Testing

In this supplement file, additional technical support is provided.

- Section S.1 provides Condition H that is used in Theorem 4.1.
- Section S.2 summarizes the remainder terms in Section 3.
- Section S.3 provides technical proofs of the lemmas in Section 6.3 in the main paper.
- Section S.4 provides technical proofs of the lemmas in Section 6.4 in the main paper.
- Section S.5 provides the proof of Theorem 4.1.
- Section S.6 provides the validity of the quadratic approximation of the PLRT test.
- Section S.7 provides a set of auxiliary results such as concentration inequalities.

S.1. *Condition H.* Condition H consists of:

(i) $h^{m-1/2} < \min\{1/(6c_K), 1/(4c_K e)\}$;

(ii) $h^{m+1/2} \leq 8c_K/5$;

(iii) $h^{1/2} r_n \leq 1$;

(iv) $72c_K^2(4e + M)r_n h^{-1/2} \leq 1$;

(v) $288c_K^4(4e + M)(Mnhr_n^2 + \log n)h^{-3/2}r_n \leq 1$;

(vi) $c_K^2(Mnhr_n^2 + \log n)\sqrt{Mh^{-1}r_n^2}A(h) + (48e^2)^{1/4}\exp(-(Mnhr_n^2 + \log n)/4)c_K h^{-1/2} \leq 1/6$,

where $A(h) = A(h, 2)$ is a known function of $h$ given in .

S.2. *Remainder terms.* In this subsection, we summarize the remainder terms appeared in Section 3.

*Simple hypothesis testing:.*

$$R_{0,n}(M) = \left(\frac{c_K^2}{\sqrt{2}\,n^{3/2}h} + \frac{c_K}{\sqrt{2}\,n^{3/2}\sqrt{h}} + \frac{6\sqrt{2}c_K^2}{n^{5/2}}\right)\sqrt{M}$$
$$+ \left(\frac{\sqrt[4]{2}\,c_K}{n^{7/4}\sqrt{h}} + \frac{8c_K}{\sqrt[4]{2}\,n^{5/4}\sqrt{h}}\right)M^{3/4} + \left(\frac{4}{n} + R_{3,n}^2(M) + \frac{8c_K^2}{n^{3/2}h} + \frac{4c_K^2}{2n^2h}\right)M,$$

$$R_{1,n}(M) = R_{0,n}(M) + 2R_{3,n}(M)\left(\left(\frac{1}{n}\sum_{\nu \geq 1}\frac{1}{1 + \lambda\rho_\nu}\right)^{1/2} + \frac{4\rho_K}{n\sqrt{h}}\sqrt{M} + R_{0,n}(M)\right).$$

$$R_{2,n}(M) = \frac{2}{n}(M^{3/4} + M) + R_{0,n}(M) + R_{4,n}(M),$$

$$R_{3,n}(M) = c_K^2\sqrt{M}n^{-1/2}h^{-1}A(h)\,\delta_n(M, (nh)^{-1/2}),$$

$$R_{4,n}(M) = R_{3,n}^2(M) + 2R_{3,n}(M)\left(1 + \left(\frac{1}{n}\sum_{\nu \geq 1}\frac{1}{1 + \lambda\rho_\nu}\right)^{1/2} + \frac{4\rho_K}{n\sqrt{h}}\sqrt{M} + R_{0,n}(M)\right).$$

*Composite hypothesis testing:.*

$$R_{0,n}^c(M) = \Big(\frac{c_K^2}{\sqrt{2}\,n^{3/2}h} + \frac{c_K}{\sqrt{2}\,n^{3/2}\sqrt{h}} + \frac{6\sqrt{2}c_K^2}{n^{5/2}}\Big)\sqrt{M}$$
$$+ \Big(\frac{\sqrt[4]{2}\,c_K}{n^{7/4}\sqrt{h}} + \frac{8c_K}{\sqrt[4]{2}\,n^{5/4}\sqrt{h}}\Big)M^{3/4} + \Big(\frac{4}{n} + \frac{8c_K^2}{n^{3/2}h} + \frac{4c_K^2}{2n^2h}\Big)M,$$

$$R_{1,n}^c(M) = R_{0,n}^c(M) + R_{4,n}^c(M) + R_{3,n}^2(M) + 2R_{3,n}(M)\Big(\Big(\frac{1}{n}\sum_{\nu\geq 1}\frac{1}{1+\lambda\rho_\nu}\Big)^{1/2}$$
$$+ \frac{4\rho_K}{n\sqrt{h}}\sqrt{M} + R_{0,n}(M) + R_{3,n}^c\Big),$$

$$R_{2,n}^c(M) = \frac{2}{n}(M^{3/4}+M) + R_{0,n}(M) + \big((R_{3,n}(M)+R_{3,n}^c(M))\big)^2$$
$$+ 2\big(R_{3,n}(M)+R_{3,n}^c(M)\big)\left(1 + \Big(\frac{1}{n}\sum_{\nu\geq 1}\frac{1}{1+\lambda\rho_\nu}\Big)^{1/2} + \frac{4\rho_K}{n\sqrt{h}}\sqrt{M} + R_{0,n}(M)\right),$$

$$R_{3,n}^c(M) = \frac{2}{\sqrt{n}}\big(1+2\sqrt{2M}+M\big)^{1/2}\Big(1-\frac{1}{n}-\sqrt{\frac{M}{n}}-\frac{M}{n}\Big)^{-1/2} + \frac{4\sqrt{2}}{\sqrt{n}}\sqrt{M}\Big(1-\frac{1}{n}-\sqrt{\frac{M}{n}}-\frac{M}{n}\Big)^{-1},$$

$$R_{4,n}^c(M) = 4\Big(\frac{1}{n}+\frac{\sqrt{2M}}{n}+\frac{M}{n}\Big)\Big(1+\Big(1-\frac{1}{n}-\sqrt{\frac{M}{n}}-\frac{M}{n}\Big)^{-1}\Big).$$

S.3. *Proof of the lemmas in Section 6.3.* We first prove some relevant large deviation inequalities, which are needed in the proofs of these lemmas.

S.3.1. *Some large deviation inequalities.* First, we consider the quadratic form

$$Q_n := \frac{1}{n^2}\sum_{i,j=1}^n \epsilon_i\epsilon_j K(X_i,X_j) = \underbrace{\frac{1}{n^2}\sum_{i=1}^n \epsilon_i^2 K(X_i,X_i)}_{V_n} + \underbrace{\frac{2}{n^2}\sum_{1\leq i<j\leq n}\epsilon_i\epsilon_j K(X_i,X_j)}_{U_n},$$

where $V_n$ collects all diagonal terms and $U_n$ collects non-diagonal terms. We provide concentration inequalities for $V_n$ and $U_n$ separately. Here, we may assume that $\epsilon_i$ are iid 1-sub-Gaussian random variables.

LEMMA S.6. *For any $x \geq 0$, it holds that*

$$P\Big(\Big|V_n - \frac{1}{n}\sum_{\nu\geq 1}\frac{1}{1+\lambda\rho_\nu}\Big| \geq \Big(\frac{2c_K^2}{nh}+\frac{c_K}{n\sqrt{h}}\Big)\sqrt{\frac{M}{2n}} + \frac{2c_K}{n\sqrt{h}}\Big(\frac{M}{2n}\Big)^{3/4} + \frac{4c_K^2}{nh}\frac{M}{2n}\Big) \leq 5\exp(-M),$$

*where $c_K = \sup_{x\in\mathbb{I}}h^{1/2}\|K_x\| = \sup_{x\in\mathbb{I}}h^{1/2}K^{1/2}(x,x)$.*

LEMMA S.7. *For any $M \geq 0$, it holds that*

$$P\Big(|U_n| \geq \Big(\frac{4\rho_K}{n\sqrt{h}}+\frac{6\sqrt{2}c_K^2}{n^{5/2}}\Big)\sqrt{M} + \frac{8c_K}{\sqrt[4]{2}\,n^{5/4}\sqrt{h}}M^{3/4} + \Big(\frac{4}{n}+\frac{8c_K^2}{n^{3/2}h}\Big)M\Big) \leq 8\exp(-M),$$

*where $\rho_K^2 = h\,E[K^2(X_1,X_2)]$.*

The following lemma provides large deviation bound for $n^{-1} \sum_{i=1}^{n} \epsilon_i f(X_i)$ over $f \in S^m(\mathbb{I})$.

LEMMA S.8. *For any $M \geq 0$, it holds that*

$$\sup_{f \in S^m(\mathbb{I})} P\left(\Big| \frac{1}{n} \sum_{i=1}^{n} \epsilon_i f(X_i) \Big| \geq \frac{\|f\|}{\sqrt{n}} \sqrt{2M} + \frac{2\|f\|^{1/2}}{n} M^{3/4} + \frac{2\|f\|_{\sup}^{1/2}}{n} M\right) \leq 2 \exp(-M).$$

S.3.2. *Proof of Lemma 6.2.* Direct calculations yields that under $H_0$, the remainder term in equation (3.8) is

$$\|\widehat{f}_{n,\lambda} - f_0 - S_{n,\lambda}(f_0)\|^2 + 2\langle \widehat{f}_{n,\lambda} - f_0 - S_{n,\lambda}(f_0), S_{n,\lambda}(0)\rangle.$$

Then the claimed bound follows by applying Lemma S.6 and Lemma S.7 to the first term in equation (3.8) and using the bound in (S.2) with $f = f_0$ for the above remainder term (see the proof of Lemma 6.3 below).

S.3.3. *Proof of Lemma 6.3.* Direct calculations gives

$$\widetilde{T}_{n,\lambda} = \|\widehat{f}_{n,\lambda} - f - S_{n,\lambda}(f) + f + S_{n,\lambda}(f) - (I - \mathcal{P}_\lambda)f_0\|^2 - \frac{1}{n} \sum_{\nu \geq 1} \frac{1}{1 + \lambda \rho_\nu}$$

$$= \underbrace{\|f + S_{n,\lambda}(f) - (I - \mathcal{P}_\lambda)f_0\|^2 - \frac{1}{n} \sum_{\nu \geq 1} \frac{1}{1 + \lambda \rho_\nu}}_{I_n}$$

$$+ \underbrace{\|\widehat{f}_{n,\lambda} - f - S_{n,\lambda}(f)\|^2 + 2\langle \widehat{f}_{n,\lambda} - f - S_{n,\lambda}(f), f - (I - \mathcal{P}_\lambda)f_0 + S_{n,\lambda}f\rangle}_{R_n}.$$

First, we provide a non-asymptotic bound for the remainder term $R_n$. By Theorem 3.4 (choose $r^2 = (nh)^{-1}$), for any $M > 0$

$$\text{(S.1)} \qquad \sup_{f \in S^m(\mathbb{I})} P_f^n\left(\|\widehat{f}_{n,\lambda} - f - S_{n,\lambda}(f)\| > R_{3,n}(M)\right) \leq 2 \exp(-M),$$

where remainder term

$$R_{3,n}(M) = c_K^2 \sqrt{M} n^{-1/2} h^{-1} A(h)\, \delta_n(M, (nh)^{-1/2}),$$

with $\delta_n(M, r) = 2h^m + c_K(\sqrt{2M}r + (nh)^{-1/2})$ the "effective" estimation convergence rate. From an asymptotic perspective, this remainder term $R_{3,n}(M)$ corresponds to the higher-order error— converging to zero faster than the estimation convergence rate $\delta_n(M, (nh)^{-1/2})$.

Combining Lemma S.6, Lemma S.7 and inequality (S.1), we obtain the following large deviation bound for the remainder term $R_n$,

$$\text{(S.2)} \qquad P\left(|R_n| \geq R_{4,n}(M)\right) \leq 15 \exp(-M),$$

where the remainder term

$$R_{4,n}(M) = R_{3,n}^2(M) + 2R_{3,n}(M)\Big(1 + \Big(\frac{1}{n}\sum_{\nu\geq1}\frac{1}{1+\lambda\rho_\nu}\Big)^{1/2} + \frac{4\rho_K}{n\sqrt{h}}\sqrt{M} + R_{0,n}(M)\Big).$$

Next, we consider the leading term $I_n$. Simple calculation suggests

$$I_n = \|(I-\mathcal{P}_\lambda)(f-f_0)\|^2 + 2\langle(I-\mathcal{P}_\lambda)(f-f_0), S_{n,\lambda}(0)\rangle + \|S_{n,\lambda}(0)\|^2 - \frac{1}{n}\sum_{\nu\geq1}\frac{1}{1+\lambda\rho_\nu}$$

$$= \|(I-\mathcal{P}_\lambda)(f-f_0)\|^2 + \underbrace{2\langle\frac{1}{n}\sum_{i=1}^n\epsilon_i K_{X_i},(I-\mathcal{P}_\lambda)(f-f_0)\rangle + \|S_{n,\lambda}(0)\|^2 - E\|S_{n,\lambda}(0)\|^2}_{W_n}.$$

The first and the second terms in the stochastic component $W_n$ ($EW_n = 0$) can be bounded by Lemma S.8 and Lemma 6.2, respectively, yielding

(S.3)
$$P\Big(|W_n| \geq \Big(\frac{2\sqrt{2}\,\|(I-\mathcal{P}_\lambda)(f-f_0)\|}{\sqrt{n}} + \frac{4\rho_K}{n\sqrt{h}}\Big)\sqrt{M} + R_{0,n}(M) + \frac{2}{n}(M^{3/4}+M)\Big) \leq 15\exp(-M),$$

where we have used the fact that $\|f-f_0\|$, $\|f-f_0\|_{\sup} \leq 1$ for $f-f_0 \in H^m(1)$.

Combining (S.2) and (S.3), we obtain

$$\sup_{f\in H^m(1)} P_f^n\Big(|\widetilde{T}_{n,\lambda} - \|f-f_0 + \mathcal{P}_\lambda f\|^2| \geq \Big(\frac{2\sqrt{2}\,\|(I-\mathcal{P}_\lambda)(f-f_0)\|}{\sqrt{n}} + \frac{4\rho_K}{n\sqrt{h}}\Big)\sqrt{M}$$

$$+ \frac{2}{n}(M^{3/4}+M) + R_{0,n}(M) + R_{4,n}(M)\Big) \leq 30\exp(-M),$$

which yields the claimed result.

S.4. *Proof of the lemmas in Section 6.4.* In this subsection, we prove the two lemmas in Section 6.4 . We use $P_X$ to denote the marginal probability measure of the covariate $\{X_i\}_{i=1}^n$.

S.4.1. *Proof of Lemma 6.4.* By Cauchy-Schwarz inequality, we have

$$\|f_\epsilon\|_2^2 = \widehat{\alpha}^2 + \widehat{\alpha}\widehat{\beta} + \frac{1}{3}\widehat{\beta}^2 \leq 2\widehat{\alpha} + 2\widehat{\beta}^2.$$

Since $(\widehat{\alpha}, \widehat{\beta})^T = (D_X^T D_X)^{-1} D_X^T \epsilon$, we can further derive that

(S.4)
$$\|f_\epsilon\|_2^2 \leq 2\epsilon^T D_X (D_X^T D_X)^{-2} D_X^T \epsilon.$$

After simple calculations, it is easy to show that the largest eigenvalue $\lambda_X$ of matrix $(D_X^T D_X)^{-1}$ is

$$\frac{n + \sum_{i=1}^n X_i^2 + \sqrt{\Big(\sum_{i=1}^n X_i^2 - n\Big)^2 + 4\Big(\sum_{i=1}^n X_i\Big)^2}}{2\Big(n\sum_{i=1}^n X_i^2 - \Big(\sum_{i=1}^n X_i\Big)^2\Big)} := \frac{N_n}{D_n}.$$

By Hanson-Wright inequality, the denominator $D_n$ satisfies

$$P_X\big(D_n \leq 2n(n-1) - 2n\sqrt{(n-1)M} - 2nM\big) \leq \exp(-M)$$

for all $M > 0$. By Cauchy-Schwarz inequality, the numerator $N_n$ satisfies

$$N_n \leq 2n + 2\sum_{i=1}^n X_i^2 \leq 4n.$$

Combining the last two displays, we obtain

$$P_X\Big(\lambda_X \geq \frac{2}{n}\Big(1 - \frac{1}{n} - \sqrt{\frac{M}{n}} - \frac{M}{n}\Big)^{-1}\Big) \leq \exp(-M)$$

for all $M \leq n/4$.

Combining the above display with inequality (S.4), we obtain that it holds with probability at least $1 - \exp(-M)$ that

$$\|f_\epsilon\|_2^2 \leq \frac{2}{n}\epsilon^T D_X (D_X^T D_X)^{-1} D_X^T \epsilon \left(1 - \frac{1}{n} - \sqrt{\frac{x}{n}} - \frac{x}{n}\right)^{-1}.$$

Moreover, since $D_X(D_X^T D_X)^{-1}D_X^T$ is a projection matrix with rank 2, we have, by Hanson-Wright inequality, that for any $M > 0$

$$P\Big(\frac{1}{n}\epsilon^T D_X (D_X^T D_X)^{-1} D_X^T \epsilon \geq \frac{2}{n} + \frac{2\sqrt{2M}}{n} + \frac{2M}{n}\Big) \leq \exp(-M).$$

Combining the last two displays, we obtain that for all $M \leq n/4$

$$P\Big\{\|f_\epsilon\|_2^2 + \frac{2}{n}\epsilon^T D_X (D_X^T D_X)^{-1} D_X^T \epsilon \geq 4\left(\frac{1}{n} + \frac{\sqrt{2M}}{n} + \frac{M}{n}\right)\Big(1 + \Big(1 - \frac{1}{n} - \sqrt{\frac{M}{n}} - \frac{M}{n}\Big)^{-1}\Big)\Big\}$$

$$\text{(S.5)} \hspace{6cm} \leq 2\exp(-M).$$

Finally, the claimed result is a direct sequence by applying Lemma S.6, Lemma S.7, equation (3.19) and equation (S.5) to the first four terms, and the bound (S.8) in the proof of Lemma 6.5 to the remainder term in decomposition (3.19).

S.4.2. *Proof of Lemma 6.5.* Write $f(x) = f^L(x) + f^\perp(x)$, where $f^L(x) = \alpha + \beta x$ is the linear part of $f$ and function $f^\perp$ satisfies $\int_0^1 f^\perp(x)\,dP(x) = \int_0^1 (f^\perp)'(x)\,dP(x) = 0$. Direct calculation gives $\widehat{f}_n^{H_0} = f_L + f_X^\perp + f_\epsilon$, and

$$\widetilde{T}_{n,\lambda}^{com} = \|\widehat{f}_{n,\lambda} - f - S_{n,\lambda}(f) + f + S_{n,\lambda}(f) - f^L - f_X^\perp - f_\epsilon\|^2 - \frac{1}{n}\sum_{\nu \geq 1}\frac{1}{1 + \lambda\rho_\nu}$$

$$= \|f^\perp + S_{n,\lambda}(0) - \mathcal{P}_\lambda f - f_X^\perp - f_\epsilon\|^2 - E\|S_{n,\lambda}(0)\|^2$$

$$+ \|\widehat{f}_{n,\lambda} - f - S_{n,\lambda}(f)\|^2 + 2\langle \widehat{f}_{n,\lambda} - f - S_{n,\lambda}(f), f + S_{n,\lambda}f - f^L - f_X^\perp - f_\epsilon\rangle,$$

where $f_X^\perp(x) = (1, x)^T (D_X^T D_X)^{-1} D_X^T f^\perp(X_1^n)$ and $f^\perp(X_1^n) = (f^\perp(X_1), \ldots, f^\perp(X_n))^T$. Use the property that $\mathcal{P}_\lambda g = 0$ for all linear functions $g$, we further obtain

$$\widetilde{T}_{n,\lambda}^{com} = \underbrace{\|f^\perp - \mathcal{P}_\lambda f^\perp + S_{n,\lambda}(0) - f_X^\perp - f_\epsilon\|^2 - E\|S_{n,\lambda}(0)\|^2}_{I_n^c}$$

$$+ \underbrace{\|\widehat{f}_{n,\lambda} - f - S_{n,\lambda}(f)\|^2 + 2\langle \widehat{f}_{n,\lambda} - f - S_{n,\lambda}(f), f^\perp - \mathcal{P}_\lambda f^\perp + S_{n,\lambda}(0) - f_X^\perp - f_\epsilon\rangle}_{R_n^c}.$$

We first provide a non-asymptotic bound for the remainder term $R_n^c$. Similar to the proof of Lemma 6.3, by Theorem 3.4 (choose $r^2 = (nh)^{-1}$), we have for any $M > 0$

$$(S.6) \qquad \sup_{f \in S^m(\mathbb{I})} P_f^n\left(\|\widehat{f}_{n,\lambda} - f - S_{n,\lambda}(f)\| > R_{3,n}(M)\right) \le 2\exp(-M).$$

According to the proof of Lemma 6.3, we have

$$(S.7) \qquad P\left\{\|f_\epsilon\|_2^2 \ge 4\left(\frac{1}{n} + \frac{\sqrt{2M}}{n} + \frac{M}{n}\right)\left(1 - \frac{1}{n} - \sqrt{\frac{M}{n}} - \frac{M}{n}\right)^{-1}\right\} \le 2\exp(-M),$$

for all $M \le n/4$. The following lemma provides non-asymptotic bounds for $\|f_X^\perp\|$, whose proof is provided in Appedix S.7.4.

LEMMA S.9.    *For any $M \in (0, n/4)$, we have*

$$P\left\{\|f_X^\perp\|^2 \ge \frac{32}{n}M\left(1 - \frac{1}{n} - \sqrt{\frac{M}{n}} - \frac{M}{n}\right)^{-2}\right\} \le 5\exp(-M).$$

Combining (S.6), (S.7), Lemma S.9 and Lemma 6.2, we obtain the following tail bound for the remainder term $R_n^c$,

$$(S.8) \qquad P\left(|R_n^c| \ge R_{5,n}^c(M)\right) \le 24\exp(-M),$$

where the remainder term

$$R_{5,n}^c(M) = R_{3,n}^2(M) + 2R_{3,n}(M)\left(\|f^\perp - \mathcal{P}_\lambda f^\perp\| + \left(\frac{1}{n}\sum_{\nu \ge 1}\frac{1}{1+\lambda\rho_\nu}\right)^{1/2} + \frac{4\rho_K}{n\sqrt{h}}\sqrt{M} + R_{0,n}(M) + R_{3,n}^c\right),$$

with

$$R_{3,n}^c = 2\left(\frac{1}{n} + \frac{\sqrt{2M}}{n} + \frac{M}{n}\right)^{1/2}\left(1 - \frac{1}{n} - \sqrt{\frac{M}{n}} - \frac{M}{n}\right)^{-1/2} + \frac{4\sqrt{2}}{\sqrt{n}}\sqrt{M}\left(1 - \frac{1}{n} - \sqrt{\frac{M}{n}} - \frac{M}{n}\right)^{-1}.$$

Next, we consider the leading term $I_n^c$. Simple calculation suggests

$$I_n = \|f^\perp - \mathcal{P}_\lambda f^\perp\|^2 + 2\langle f^\perp - \mathcal{P}_\lambda f^\perp, S_{n,\lambda}(0)\rangle + \|S_{n,\lambda}(0)\|^2 - E\|S_{n,\lambda}(0)\|^2$$

$$+ \|f_X^\perp + f_\epsilon\|^2 - 2\langle f_X^\perp + f_\epsilon, f^\perp - \mathcal{P}_\lambda f^\perp + S_{n,\lambda}(0)\rangle$$

$$= \|f^\perp - \mathcal{P}_\lambda f^\perp\|^2 + \underbrace{2\langle\frac{1}{n}\sum_{i=1}^n \epsilon_i K_{X_i}, f^\perp - \mathcal{P}_\lambda f^\perp\rangle + \|S_{n,\lambda}(0)\|^2 - E\|S_{n,\lambda}(0)\|^2}_{W_{1,n}^c}$$

$$+ \underbrace{\|f_X^\perp + f_\epsilon\|^2 - 2\langle f_X^\perp + f_\epsilon, f^\perp - \mathcal{P}_\lambda f^\perp + S_{n,\lambda}(0)\rangle}_{W_{2,n}^c}.$$

Similar to the proof of Lemma 6.3, we have the following tail bound for $W_{1,n}^c$,

(S.9)
$$P\Big(|W_{1,n}^c| \geq \Big(\frac{2\sqrt{2}\,\|f^\perp - \mathcal{P}_\lambda f^\perp\|}{\sqrt{n}} + \frac{4\rho_K}{n\sqrt{h}}\Big)\sqrt{M} + R_{0,n}(M) + \frac{2}{n}(M^{3/4} + M)\Big) \leq 15\exp(-M).$$

Combining (S.7), Lemma S.9 and Lemma 6.2, we obtain the following tail bound for $W_{2,n}^c$,

(S.10)
$$P\Big(|W_{2,n}^c| \geq R_{5,n}^2 + 2R_{3,n}^c\,\Big(\|f^\perp - \mathcal{P}_\lambda f^\perp\| + \Big(\frac{1}{n}\sum_{\nu \geq 1}\frac{1}{1 + \lambda\rho_\nu}\Big)^{1/2} + \frac{4\rho_K}{n\sqrt{h}}\sqrt{M} + R_{0,n}(M)\Big)\Big) \leq 20\exp(-M).$$

Combining (S.2) and (S.3), we obtain that for all $x \in (0, n/4]$,

$$\begin{aligned}
\sup_{f \in H^m(1)} P_f^n\Big(|\widetilde{T}_{n,\lambda}^{com} - \|f^\perp - \mathcal{P}_\lambda f^\perp\|^2| \geq &\ \Big(\frac{2\sqrt{2}\,\|f^\perp - \mathcal{P}_\lambda f^\perp\|}{\sqrt{n}} + \frac{4\rho_K}{n\sqrt{h}}\Big)\sqrt{M} \\
&+ R_{5,n}^2 + 2R_{3,n}^c\,\Big(\|f^\perp - \mathcal{P}_\lambda f^\perp\| + \Big(\frac{1}{n}\sum_{\nu \geq 1}\frac{1}{1 + \lambda\rho_\nu}\Big)^{1/2} + \frac{4\rho_K}{n\sqrt{h}}\sqrt{M} + R_{0,n}(M)\Big) \\
&+ \frac{2}{n}(M^{3/4} + M) + R_{0,n}(M) + R_{4,n}^c(M)\Big) \leq 60\exp(-M),
\end{aligned}$$

which yields the claimed result.

S.5. *Proof of Theorem 4.1.* The proof of (i) is finished in two parts.

**Part I**: Let $f \in H^m(1)$ be the parameter based on which the data are drawn. Define an operator mapping $S^m(\mathbb{I})$ to $S^m(\mathbb{I})$:

$$T_{1f}(g) = g + S_\lambda(f + g),\ \ g \in S^m(\mathbb{I}).$$

First observe that

$$\|S_\lambda(f)\| = \|\mathcal{P}_\lambda f\| = \sup_{\|g\| = 1}|\langle \mathcal{P}_\lambda f, g \rangle| \leq \sqrt{\lambda b^2} = bh^m.$$

Let $r_{1n} = 2h^m$. Let $\mathbb{B}(r_{1n}) = \{g \in S^m(\mathbb{I}) : \|g\| \leq r_{1n}\}$ be the $r_{1n}$-ball. For any $g \in \mathbb{B}(r_{1n})$, using $DS_\lambda(f) = -id$ (see [20]) and $\|g\|_{\sup} \leq c_K h^{-1/2} r_{1n} = 2c_K h^{m-1/2} \leq 1$, it is easy to see that

$$\begin{aligned}
\|T_{1f}(g)\| &\leq \|g + S_\lambda(f + g) - S_\lambda(f)\| + \|S_\lambda(f)\| \\
&= \|g + DS_\lambda(f)g + \int_0^1\int_0^1 sD^2S_\lambda(f + ss'g)gg\,ds\,ds'\| + \|S_\lambda(f)\| \\
&= \|\int_0^1\int_0^1 sD^2S_\lambda(f + ss'g)gg\,ds\,ds'\| + \|S_\lambda(f)\| \\
&= \|\int_0^1\int_0^1 sE\{\ell^{(3)}(Y - f(X) - ss'g(X))g(X)^2K_X\}ds\,ds'\| + r_{1n}/2 \\
&\leq c_K h^{-1/2}\int_0^1\int_0^1 sE\{\sup_{|a| \leq 1}|\ell^{(3)}(\eta - a)| \cdot g(X)^2\} + r_{1n}/2 \\
&\leq c_K C_0 C_1 \sigma^{-2} h^{-1/2}\|g\|^2/2 + r_{1n}/2 \\
&\leq c_K C_0 C_1 \sigma^{-2} h^{-1/2} r_{1n}^2/2 + r_{1n}/2 \\
&= c_K C_0 C_1 \sigma^{-2} h^{m-1/2} r_{1n} + r_{1n}/2 \leq 3r_{1n}/4,
\end{aligned}$$

where the last step follows from the assumption $c_K C_0 C_1 \sigma^{-2} h^{m-1/2} < 1/4$. Therefore, $T_{1f}$ maps $\mathbb{B}(r_{1n})$ to itself.

For any $g_1, g_2 \in \mathbb{B}(r_{1n})$, denote $g = g_1 - g_2$. Note that for any $0 \le s \le 1$, $\|g_2 + sg\|_{\sup} \le c_K h^{-1/2}\|g_2 + sg\| \le 2c_K h^{-1/2} r_{1n} = 4c_K h^{-1/2} h^m \le 2/3$. According to [20], we have that for any $f \in S^m(\mathbb{I})$, $\|f\|_{\sup} \le c_K h^{-1/2}\|f\|$. Therefore,

$$
\begin{aligned}
&\|T_{1f}(g_1) - T_{1f}(g_2)\| \\
=\ & \|g_1 - g_2 + S_\lambda(f + g_1) - S_\lambda(f + g_2)\| \\
=\ & \left\|g_1 - g_2 + \int_0^1 DS_\lambda(f + g_2 + sg)g\,ds\right\| \\
=\ & \left\|\int_0^1 [DS_\lambda(f + g_2 + sg) - DS_\lambda(f)]g\,ds\right\| \\
=\ & \left\|\int_0^1 \int_0^1 D^2 S_\lambda(f + s'(g_2 + sg))(g_2 + sg)g\,ds\,ds'\right\| \\
\le\ & \int_0^1 \int_0^1 \|E\{\ell^{(3)}(Y - f(X) - s'(g_2(X) + sg(X)))(g_2(X) + sg(X))g(X)K_X\}\| \\
\le\ & c_K h^{-1/2} \int_0^1 E\{\sup_{|a|\le 1}|\ell^{(3)}(\eta - a)| \cdot |g_2(X) + sg(X)| \cdot |g(X)|\} \\
\le\ & c_K h^{-1/2} \sigma^{-2} C_0 C_1 \int_0^1 \|g_2 + sg\|\,ds\|g\| \\
\le\ & 2c_K \sigma^{-2} C_0 C_1 h^{-1/2} r_{1n}\|g\| \\
=\ & 4c_K C_0 C_1 \sigma^{-2} h^{m-1/2}\|g\| = 4c_K C_0 C_1 \sigma^{-2} h^{m-1/2}\|g_1 - g_2\|.
\end{aligned}
$$

Since $4c_K C_0 C_1 \sigma^{-2} h^{m-1/2} < 1$, this shows that $T_{1f}$ is a contraction mapping which maps $\mathbb{B}(r_{1n})$ into $\mathbb{B}(r_{1n})$. By contraction mapping theorem (see [19]), $T_{1f}$ has a unique fixed point $g' \in \mathbb{B}(r_{1n})$ satisfying $T_{1f}(g') = g'$. Let $f_\lambda = f + g'$. Then $S_\lambda(f_\lambda) = 0$ and $\|f_\lambda - f\| \le r_{1n}$.

**Part II**: For any $f \in H^m(1)$, under (2.1) with $f$ being the truth, let $f_\lambda$ be the function obtained in **Part I**. Obviously, $\|f_\lambda - f\|_{\sup} \le c_K h^{-1/2} r_{1n} = 2c_K h^{m-1/2} \le 1/3$. Then it can be shown that for all $g_1, g_2 \in S^m(\mathbb{I})$,

$$
\begin{aligned}
&|[DS_\lambda(f_\lambda) - DS_\lambda(f)]g_1 g_2| \\
=\ & |E\{(\ell^{(2)}(Y - f_\lambda(X)) - \ell^{(2)}(Y - f(X)))g_1(X)g_2(X)\}| \\
\le\ & |E\{\sup_{|a|\le 1}|\ell^{(3)}(\eta + a)| \cdot |f_\lambda(X) - f(X)| \cdot |g_1(X)g_2(X)|\}| \\
\le\ & C_0 C_1 c_K h^{-1/2} r_{1n} \sigma^{-2}\|g_1\| \cdot \|g_2\| \le \|g_1\| \cdot \|g_2\|/2,
\end{aligned}
$$

where the last inequality follows by $C_0 C_1 c_K h^{-1/2} r_{1n} \sigma^{-2} = 2C_0 C_1 c_K \sigma^{-2} h^{m-1/2} \le 1/2$. Together with the fact $DS_\lambda(f) = -id$, we get that the operator norm $\|DS_\lambda(f_\lambda) + id\|_{\text{operator}} \le 1/2$. This implies that $DS_\lambda(f_\lambda)$ is invertible with operator norm within $[1/2, 3/2]$, and hence, $\|DS_\lambda(f_\lambda)^{-1}\|_{\text{operator}} \le 2$.

Define an operator $T_{2f}(g) = g - [DS_\lambda(f_\lambda)]^{-1}S_{n,\lambda}(f_\lambda + g), \ g \in S^m(\mathbb{I})$ and rewrite it as

$$
\begin{aligned}
T_{2f}(g) &= -DS_\lambda(f_\lambda)^{-1}[DS_{n,\lambda}(f_\lambda)g - DS_\lambda(f_\lambda)g] \\
&\quad -DS_\lambda(f_\lambda)^{-1}[S_{n,\lambda}(f_\lambda + g) - S_{n,\lambda}(f_\lambda) - DS_{n,\lambda}(f_\lambda)g] \\
&\quad -DS_\lambda(f_\lambda)^{-1}S_{n,\lambda}(f_\lambda).
\end{aligned}
$$

Denote the above three terms by $I_{1f}, I_{2f}, I_{3f}$, respectively.

For $i = 1, \ldots, n$, let $R_i = \ell^{(1)}(Y_i - f_\lambda(X_i))K_{X_i} - E_f\{\ell^{(1)}(Y - f_\lambda(X))K_X\}$. By direct calculations, it can be shown that

$$
\begin{aligned}
&\|E_f\{\ell^{(1)}(Y - f_\lambda(X))K_X\}\| \\
={}& \sup_{\|g\|=1} |\langle E_f\{\ell^{(1)}(Y - f_\lambda(X))K_X\}, g\rangle| \\
={}& \sup_{\|g\|=1} |E_f\{\ell^{(1)}(Y - f_\lambda(X))g(X)\}| \\
={}& \sup_{\|g\|=1} \left| E_f\left\{(\ell^{(1)}(\eta) + \ell^{(2)}(\eta)(f(X) - f_\lambda(X)) + \frac{1}{2}\ell^{(3)}(\eta + s(f(X) - f_\lambda(X)))(f(X) - f_\lambda(X))^2)g(X)\right\} \right| \\
={}& \sup_{\|g\|=1} \left| E_f\left\{\ell^{(2)}(\eta)(f(X) - f_\lambda(X))g(X)\right\} + \frac{1}{2}E_f\left\{\ell^{(3)}(\eta + s(f(X) - f_\lambda(X)))(f(X) - f_\lambda(X))^2g(X)\right\} \right| \\
={}& \sup_{\|g\|=1} \left| -\sigma^2 E_f\left\{(f(X) - f_\lambda(X))g(X)\right\} + \frac{1}{2}E_f\left\{\ell^{(3)}(\eta + s(f(X) - f_\lambda(X)))(f(X) - f_\lambda(X))^2g(X)\right\} \right| \\
\leq{}& \|f_\lambda - f\| + \frac{1}{2}c_K h^{-1/2}E_f\{\sup_{|a|\leq 1}|\ell^{(3)}(\eta + a)| \cdot (f_\lambda(X) - f(X))^2\} \\
\leq{}& r_{1n} + \frac{1}{2}c_K h^{-1/2}C_0C_1\sigma^{-2}r_{1n}^2 \leq 5r_{1n}/4,
\end{aligned}
$$

where the second last inequality follows by Assumption A2, i.e., $E\{\sup_{|a|\leq 1}|\ell^{(3)}(\eta + a)|\} \leq C_0C_1$, and the last inequality follows by condition $h^{m-1/2} < \sigma^2/(4c_K C_0 C_1)$. Therefore, it can be shown by $\|f_\lambda - f\|_{\sup} \leq 1/3$ that

$$
\begin{aligned}
\|R_i\| &\leq \|\ell^{(1)}(Y_i - f_\lambda(X_i))K_{X_i}\| + 5r_{1n}/4 \\
&\leq \left(|\ell^{(1)}(\eta)| + \frac{1}{3}|\ell^{(2)}(\eta)| + \frac{1}{18}\sup_{|a|\leq 1}|\ell^{(3)}(\eta + a)|\right)c_K h^{-1/2} + 5r_{1n}/4 \\
&\leq c_K h^{-1/2}|\ell^{(1)}(\eta)| + \frac{1}{3}|\ell^{(2)}(\eta)| + \frac{1}{18}\sup_{|a|\leq 1}|\ell^{(3)}(\eta + a)| + 5r_{1n}/4.
\end{aligned}
$$

Using Cauchy-Schwartz inequality,

$$
E\left\{\exp\left(\frac{\|R_i\|}{2C_0c_K h^{-1/2}}\right)\right\} \leq C_1 \exp\left(\frac{5h^{m+1/2}}{4C_0c_K}\right) \leq C_1 \exp(2).
$$

Let $\delta = hr_n/(4C_0c_K)$. Recall that $h^{1/2}r_n \leq 1$ which implies $\delta \leq (4C_0c_K h^{-1/2})^{-1}$. Therefore, $E\{\exp(\delta\|R_i\|)\} \leq C_1 \exp(2)$. Moreover, for $x \geq 0$ and any constant $c > 0$, $\exp(M) - $

$1 - M \leq M^2 \exp(M)$ and $M^{-2} \exp(cx) \geq c^2 \exp(2)/4$. Let $c = (2C_0 c_K h^{-1/2})^{-1} - \delta$. Clearly, $c \geq (4C_0 c_K h^{-1/2})^{-1}$. So, we have

$$
\begin{aligned}
& E_f\{\exp(\delta\|R_i\|) - 1 - \delta\|R_i\|\} \\
\leq\ & E_f\{(\delta\|R_i\|)^2 \exp(\delta\|R_i\|)\} \\
\leq\ & \delta^2 \frac{4\exp(-2)}{c^2} E_f\{\exp(c\|R_i\|)\exp(\delta\|R_i\|)\} \\
\leq\ & \delta^2 4\exp(-2)(16C_0^2 c_K^2 h^{-1})C_1\exp(2) \leq 64 C_0^2 C_1 c_K^2 h^{-1}\delta^2.
\end{aligned}
$$

It follows by Theorem 3.2 of [18] that, for $L(M) \equiv 4C_0 c_K(4C_1 + M)$,

$$
P_f^n\left(\|\sum_{i=1}^n R_i\| \geq L(M)nr_n\right) \leq 2\exp\left(-L(M)\delta nr_n + 64C_0^2 C_1 c_K^2 nh^{-1}\delta^2\right) = 2\exp(-Mnhr_n^2).
$$

We note that the right hand side in the above inequality does not depend on $f$. It is easy to see that $S_{n,\lambda}(f_\lambda) = S_{n,\lambda}(f_\lambda) - S_\lambda(f_\lambda) = \frac{1}{n}\sum_{i=1}^n R_i$. Let

$$
\mathcal{E}_{n,1} = \{\|S_{n,\lambda}(f_\lambda)\| \leq L(M)r_n\},
$$

then $\sup_{f \in H^m(1)} P_f^n(\mathcal{E}_{n,1}^c) \leq 2\exp(-Mnhr_n^2)$.

It follows by Assumption A2 that $\sup_{f \in H^m(1)} P_f^n(\mathcal{E}_{n,2}^c) \leq 2C_1\exp(-Mnhr_n^2)$, where $\mathcal{E}_{n,2} = \cap_{i=1}^n A_i$, and $A_i = \left\{\sup_{|a|\leq 1}|\ell^{(j)}(\eta_i + a)| \leq C_0(Mnhr_n^2 + \log n),\ j = 2,3\right\}$. Define

$$
\psi_{n,f}^{(1)}(Z_i; g) = \frac{\ell^{(2)}(Y_i - f_\lambda(X_i))}{C_0(Mnhr_n^2 + \log n)} I_{A_i} c_K^{-1} h^{1/2} g(X_i)
$$

and $Z_{n,f}^{(1)}(g) = \frac{1}{\sqrt{n}}\sum_{i=1}^n [\psi_{n,f}^{(1)}(Z_i; g)K_{X_i} - E_f\{\psi_{n,f}^{(1)}(Z_i; g)K_{X_i}\}]$. It follows by Lemma 6.1 that $\sup_{f \in H^m(1)} P_f^n(\mathcal{E}_{n,3}^c) \leq 2\exp(-Mnhr_n^2)$, where $\mathcal{E}_{n,3} = \{\sup_{g \in \mathcal{G}}\|Z_{n,f}^{(1)}(g)\| \leq \sqrt{Mnhr_n^2}A(h)\}$.

For any $g \in S^m(\mathbb{I})\backslash\{0\}$, let $\bar{g} = g/d_n'$, where $d_n' = c_K h^{-1/2}\|g\|$. It is easy to see that

(S.11)                    $\|\bar{g}\|_{\sup} \leq c_K h^{-1/2}\|\bar{g}\| = c_K h^{-1/2}\|g\|/d_n' = 1,$ and

(S.12)                    $J(\bar{g}, \bar{g}) = d_n'^{-2}J(g,g) = h^{-2m}\frac{\lambda J(g,g)}{c_K^2 h^{-1}\|g\|^2} \leq c_K^{-2}h^{-2m+1}.$

Therefore, $\bar{g} \in \mathcal{G}$. Consequently, on $\mathcal{E}_{n,3}$, for any $g \in S^m(\mathbb{I})\backslash\{0\}$, we get $\|Z_{n,f}^{(1)}(\bar{g})\| \leq \sqrt{Mnhr_n^2}A(h)$, which leads to that

$$
\begin{aligned}
& \frac{1}{n}\|\sum_{i=1}^n[\ell^{(2)}(Y_i - f_\lambda(X_i))g(X_i)K_{X_i}I_{A_i} - E_f\{\ell^{(2)}(Y_i - f_\lambda(X_i))g(X_i)K_{X_i}I_{A_i}\}]\| \\
\leq\ & c_K^2 C_0(Mnhr_n^2 + \log n)\sqrt{Mh^{-1}r_n^2}A(h)\|g\|.
\end{aligned}
$$

Note that the above inequality also holds for $g = 0$.

On the other hand, for any $f, g \in S^m(\mathbb{I})$, by Cauchy-Schwartz inequality,

$$
\begin{aligned}
&\|E_f\{\ell^{(2)}(Y_i - f_\lambda(X_i))g(X_i)K_{X_i}I_{A_i^c}\}\| \\
\leq\ & E_f\{\sup_{|a| \leq 1} |\ell^{(2)}(\eta_i + a)| \cdot |g(X_i)|I_{A_i^c}\}c_K h^{-1/2} \\
\leq\ & E_f\{(\sup_{|a| \leq 1} |\ell^{(2)}(\eta_i + a)|)^2 I_{A_i^c}\}^{1/2} E\{g(X)^2\}^{1/2} c_K h^{-1/2} \\
\leq\ & E_f\{(\sup_{|a| \leq 1} |\ell^{(2)}(\eta_i + a)|)^4\}^{1/4} P_f^n(A_i^c)^{1/4} E\{g(X)^2\}^{1/2} c_K h^{-1/2} \\
\leq\ & \sigma^{-1}(48C_0^4 C_1^2)^{1/4} \exp(-(Mnhr_n^2 + \log n)/4)c_K h^{-1/2}\|g\|,
\end{aligned}
$$

where in the last inequality we have used the fact

$$
\text{(S.13)} \qquad\qquad P_f^n(A_i^c) \leq 2C_1 \exp(-(Mnhr_n^2 + \log n)).
$$

In summary, we have shown that for any $f \in H^m(1)$, on $\mathcal{E}_{n,2} \cap \mathcal{E}_{n,3}$, uniformly for $g \in S^m(\mathbb{I})$,

$$
\begin{aligned}
&\|DS_{n,\lambda}(f_\lambda)g - DS_\lambda(f_\lambda)g\| \\
=\ & \frac{1}{n}\|\sum_{i=1}^n [\ell^{(2)}(Y_i - f_\lambda(X_i))g(X_i)K_{X_i} - E_f\{\ell^{(2)}(Y_i - f_\lambda(X_i))g(X_i)K_{X_i}\}]\| \\
=\ & \frac{1}{n}\|\sum_{i=1}^n [\ell^{(2)}(Y_i - f_\lambda(X_i))g(X_i)K_{X_i}I_{A_i} - E_f\{\ell^{(2)}(Y_i - f_\lambda(X_i))g(X_i)K_{X_i}\}]\| \\
\leq\ & \frac{1}{n}\|\sum_{i=1}^n [\ell^{(2)}(Y_i - f_\lambda(X_i))g(X_i)K_{X_i}I_{A_i} - E_f\{\ell^{(2)}(Y_i - f_\lambda(X_i))g(X_i)K_{X_i}I_{A_i}\}]\| \\
& + \|E_f\{\ell^{(2)}(Y_i - f_\lambda(X_i))g(X_i)K_{X_i}I_{A_i^c}\}\| \\
\leq\ & (c_K^2 C_0(Mnhr_n^2 + \log n)\sqrt{Mh^{-1}r_n^2}A(h) \\
\text{(S.14)} \quad & + \sigma^{-1}(48C_0^4 C_1^2)^{1/4}\exp(-(Mnhr_n^2 + \log n)/4)c_K h^{-1/2})\|g\| \leq \|g\|/6.
\end{aligned}
$$

Define $T_{3f}(g) = S_{n,\lambda}(f_\lambda + g) - S_{n,\lambda}(f_\lambda) - DS_{n,\lambda}(f_\lambda)g$. Let $r_{2n} = 6L(M)r_n$. For any $g_1, g_2 \in \mathbb{B}(r_{2n})$, and $s \in \mathbb{I}$, let $g = g_1 - g_2$, then $\|g_2 + sg\|_{\sup} \leq \|g_1\|_{\sup} + \|g_2\|_{\sup} \leq 2c_K h^{-1/2}r_{2n} = 48C_0 c_K^2 (4C_1 + M)r_n h^{-1/2} < 2/3$. On $\mathcal{E}_{n,2} \cap \mathcal{E}_{n,3}$, for any $g_1, g_2 \in \mathbb{B}(r_{2n})$ and letting $g = g_1 - g_2$,

we have

$$\|T_{3f}(g_1) - T_{3f}(g_2)\|$$

$$= \|S_{n,\lambda}(f_\lambda + g_1) - S_{n,\lambda}(f_\lambda + g_2) - DS_{n,\lambda}(f_\lambda)g\|$$

$$= \|\int_0^1 \int_0^1 D^2 S_{n,\lambda}(f_\lambda + s'(g_2 + sg))(g_2 + sg)gdsds'\|$$

$$\leq \int_0^1 \int_0^1 \|D^2 S_{n,\lambda}(f_\lambda + s'(g_2 + sg))(g_2 + sg)g\|dsds'$$

$$\leq \int_0^1 \int_0^1 \|\frac{1}{n}\sum_{i=1}^n \ell^{(3)}(Y_i - f_\lambda(X_i) - s'(g_2(X_i) + sg(X_i)))(g_2(X_i) + sg(X_i))g(X_i)K_{X_i}\|dsds'$$

$$\leq \int_0^1 \int_0^1 \frac{1}{n}\sum_{i=1}^n \sup_{|a|\leq 1} |\ell^{(3)}(\eta_i + a)| \cdot \|g_2 + sg\|_{\sup} \cdot \|g\|_{\sup}\|K_{X_i}\|dsds'$$

$$\leq 2C_0(Mnhr_n^2 + \log n)(c_K h^{-1/2})^3 r_{2n}\|g\|$$

$$\text{(S.15)} \quad 48C_0^2 c_K^4 (4C_1 + M)(Mnhr_n^2 + \log n)h^{-3/2}r_n\|g\| \leq \|g_1 - g_2\|/6.$$

Taking $g_2 = 0$ in (S.15) we get that $\|T_{3f}(g_1)\| \leq \|g_1\|/6$ for any $g_1 \in \mathbb{B}(r_{2n})$. Therefore, it follows by (S.14) and (S.15) that, for any $f \in H^m(1)$, on $\mathcal{E}_n \equiv \mathcal{E}_{n,1} \cap \mathcal{E}_{n,2} \cap \mathcal{E}_{n,3}$ and for any $g \in \mathbb{B}(r_{2n})$,

$$\|T_{2f}(g)\| \leq 2(\|g\|/6 + \|g\|/6 + r_{2n}/6) \leq 2(r_{2n}/6 + r_{2n}/6 + r_{2n}/6) = r_{2n},$$

meanwhile, for any $g_1, g_2 \in \mathbb{B}(r_{2n})$, replacing $g$ by $g_1 - g_2$ in (S.14) we get that

$$\|T_{2f}(g_1) - T_{2f}(g_2)\| \leq 2(\|g_1 - g_2\|/6 + \|g_1 - g_2\|/6) = 2\|g_1 - g_2\|/3.$$

Therefore, for any $f \in H^m(1)$, on $\mathcal{E}_n$, $T_{2f}$ is a contraction mapping from $\mathbb{B}(r_{2n})$ to itself. By contraction mapping theorem, there exists uniquely an element $g'' \in \mathbb{B}(r_{2n})$ s.t. $T_{2f}(g'') = g''$. Let $\widehat{f}_{n,\lambda} = f_\lambda + g''$. Then on $\mathcal{E}_n$, $\|\widehat{f}_{n,\lambda} - f\| \leq \|f_\lambda - f\| + \|\widehat{f}_{n,\lambda} - f_\lambda\| \leq r_{1n} + r_{2n} = 2bh^m + 6L(M)r_n$. The desired conclusion follows by the trivial fact: $\sup_{f \in H^m(1)} P_f^n(\mathcal{E}_n^c) \leq (2C_1 + 4)\exp(-Mnhr_n^2)$. Proof of (i) is completed.

Next we show (ii).

For any $f \in H^m(1)$, let $\widehat{f}_{n,\lambda}$ be the penalized MLE of $f$. Let $g_n = \widehat{f}_{n,\lambda} - f$, $\delta'_n(M) = 2h^m + 6L(M)r_n$, $d'_n = c_K h^{-1/2}\delta'_n(M)$, and for $g \in \mathcal{G}$ define

$$\psi_{n,f}^{(2)}(Z_i; g) = \frac{\ell^{(1)}(\eta_i - d'_n g(X_i)) - \ell^{(1)}(\eta_i)}{C_0(Mnhr_n^2 + \log n)(c_K h^{-1/2})^2 \delta'_n(M)} I_{A_i},$$

where $A_i$ is the event defined in (i). Under the imposed conditions, we get that

$$\text{(S.16)} \quad c_K h^{-1/2}\delta'_n(M) = 2c_K h^{m-1/2} + 6L(M)c_K h^{-1/2}r_n \leq 2/3.$$

Then for any $g_1, g_2 \in \mathcal{G}$,

$$
\begin{aligned}
&|\psi_{n,f}^{(2)}(Z_i; g_1) - \psi_{n,f}^{(2)}(Z_i; g_2)| \\
=\; & \frac{1}{C_0(Mnhr_n^2 + \log n)(c_K h^{-1/2})^2 \delta_n'(M)} |\ell^{(1)}(\eta_i - d_n' g_1(X_i)) - \ell^{(1)}(\eta_i - d_n' g_2(X_i))| I_{A_i} \\
=\; & \frac{1}{C_0(Mnhr_n^2 + \log n)(c_K h^{-1/2})^2 \delta_n'(M)} \\
& \cdot |\int_0^1 \ell^{(2)}(\eta_i - d_n' g_2(X_i) + s d_n(g_2(X_i) - g_1(X_i))) d_n'(g_2(X_i) - g_1(X_i)) ds| I_{A_i} \\
\leq\; & \frac{1}{C_0(Mnhr_n^2 + \log n)(c_K h^{-1/2})^2 \delta_n'(M)} \sup_{|a| \leq 1} |\ell^{(2)}(\eta_i + a)| I_{A_i} d_n' \|g_1 - g_2\|_{\sup} \\
=\; & c_K^{-1} h^{1/2} \|g_1 - g_2\|_{\sup}.
\end{aligned}
$$

Let $\mathcal{E}_{n,4} = \{\sup_{g \in \mathcal{G}} \|Z_{n,f}^{(2)}(g)\| \leq \sqrt{Mnhr_n^2} A(h)\}$, where $Z_{n,f}^{(2)}(g) = \frac{1}{\sqrt{n}} \sum_{i=1}^n [\psi_{n,f}^{(2)}(Z_i; g) K_{X_i} - E_f^Z \{\psi_{n,f}^{(2)}(Z_i; g) K_{X_i}\}]$, $E_f^Z$ denotes the expectation with respect to $Z$ (under $P_f^n$). It follows by Lemma 6.1 that $\sup_{f \in H^m(1)} P_f^n(\mathcal{E}_{n,4}^c) \leq 2\exp(-Mnhr_n^2)$.

On the other hand, for any $g \in \mathcal{G}$, using $P_f^n(A_i^c) \leq 2C_1 \exp(-(Mnhr_n^2 + \log n))$ we get

$$
\begin{aligned}
& \|E_f^Z \{(\ell^{(1)}(\eta_i - d_n' g(X_i)) - \ell^{(1)}(\eta_i)) K_{X_i} I_{A_i^c}\}\| \\
\leq\; & E_f^Z \{\sup_{|a| \leq 1} |\ell^{(2)}(\eta_i + a)| \cdot d_n' |g(X_i)| \cdot \|K_{X_i}\| I_{A_i^c}\} \\
\leq\; & d_n c_K h^{-1/2} E_f^Z \{\sup_{|a| \leq 1} |\ell^{(2)}(\eta_i + a)| \cdot I_{A_i^c}\} \\
\leq\; & (c_K h^{-1/2})^2 \delta_n'(M) \sqrt{2C_0^2 C_1} \sqrt{2C_1} \exp(-(Mnhr_n^2 + \log n)/2) \\
=\; & 2c_K^2 C_0 C_1 h^{-1} \delta_n'(M) \exp(-(Mnhr_n^2 + \log n)/2).
\end{aligned}
$$

On $\widetilde{\mathcal{E}}_n \equiv \mathcal{E}_n \cap \mathcal{E}_{n,4}$, we have $\|g_n\| \leq \delta_n'(M)$. Let $\bar{g} = g_n/d_n'$. Then we get that

$$
\begin{aligned}
& \|S_{n,\lambda}(f + g_n) - S_{n,\lambda}(f) - (S_\lambda(f + g_n) - S_\lambda(f))\| \\
=\; & \frac{1}{n} \|\sum_{i=1}^n [(\ell^{(1)}(\eta_i - g_n(X_i)) - \ell^{(1)}(\eta_i)) K_{X_i} - E_f \{(\ell^{(1)}(\eta_i - g_n(X_i)) - \ell^{(1)}(\eta_i)) K_{X_i}\}]\| \\
\leq\; & \frac{1}{n} C_0(Mnhr_n^2 + \log n)(c_K h^{-1/2})^2 \delta_n'(M) \|\sum_{i=1}^n [\psi_{n,f}^{(2)}(Z_i; \bar{g}) K_{X_i} - E_f^Z \{\psi_{n,f}^{(2)}(Z_i; \bar{g}) K_{X_i}\}]\| \\
& + C_0(Mnhr_n^2 + \log n)(c_K h^{-1/2})^2 \delta_n'(M) \|E_f^Z \{(\ell^{(1)}(\eta_i - d_n' \bar{g}(X_i)) - \ell^{(1)}(\eta_i)) K_{X_i} I_{A_i^c}\}\| \\
\leq\; & \frac{1}{n} C_0(Mnhr_n^2 + \log n)(c_K h^{-1/2})^2 \delta_n'(M) \cdot \sqrt{n} \sqrt{Mnhr_n^2} A(h) \\
& + C_0(Mnhr_n^2 + \log n)(c_K h^{-1/2})^2 \delta_n'(M) \cdot 2c_K^2 C_0 C_1 h^{-1} \delta_n'(M) \exp(-(Mnhr_n^2 + \log n)/2) \\
=\; & c_K^2 C_0 \sqrt{M} (Mnhr_n^2 + \log n) h^{-1/2} \delta_n'(M) r_n A(h) \\
& + 2c_K^4 C_0^2 C_1 h^{-2} \delta_n'(M)^2 (Mnhr_n^2 + \log n) \exp(-(Mnhr_n^2 + \log n)/2) = \alpha_n. \quad \text{(S.17)}
\end{aligned}
$$

It is easy to show that

$$\| \int_0^1 \int_0^1 s D^2 S_\lambda(f + ss'g_n) g_n g_n ds ds' \|$$

$$= \| \int_0^1 \int_0^1 s E_f^Z \{ \ell^{(3)}(Y - f(X) - ss'g_n(X)) g_n(X)^2 K_X \} ds ds' \|$$

$$\leq c_K h^{-1/2} \int_0^1 \int_0^1 s E_f^Z \{ |\ell^{(3)}(\eta - ss'g_n(X))| g_n(X)^2 \} ds ds'$$

$$\leq \frac{1}{2} c_K h^{-1/2} E \{ \sup_{|a| \leq 1} |\ell^{(3)}(\eta + a)| g_n(X)^2 \}$$

$$\text{(S.18)} \qquad \leq \frac{1}{2} c_K h^{-1/2} C_0 C_1 \sigma^{-2} \| g_n \|^2 \leq \frac{1}{2} c_K C_0 C_1 \sigma^{-2} h^{-1/2} \delta_n'(M)^2 = \beta_n.$$

Since $S_{n,\lambda}(f + g_n) = 0$ and $DS_\lambda(f) = -id$, from (S.17) and (S.18) we have on $\widetilde{\mathcal{E}}_n$,

$$\alpha_n \geq \| S_{n,\lambda}(f) + DS_\lambda(f) g_n + \int_0^1 \int_0^1 s D^2 S_\lambda(f + ss'g_n) g_n g_n ds ds' \|$$

$$= \| S_{n,\lambda}(f) - g_n + \int_0^1 \int_0^1 s D^2 S_\lambda(f + ss'g_n) g_n g_n ds ds' \|$$

$$\geq \| S_{n,\lambda}(f) - g_n \| - \| \int_0^1 \int_0^1 s D^2 S_\lambda(f + ss'g_n) g_n g_n ds ds' \|,$$

which implies that

$$\| \widehat{f}_{n,\lambda} - f - S_{n,\lambda}(f) \| \leq c_n(M),$$

where $c_n(M) := \alpha_n + \beta_n$. Since $\sup_{f \in H^m(1)} P_f^n(\widetilde{\mathcal{E}}_n^c) \leq (2C_1 + 6) \exp(-Mnhr_n^2)$, the proof is completed.

S.6. *A quadratic approximation of the PLRT test.* For arbitrary $f \in H^m(1)$, let $\widehat{f}_{n,\lambda}$ be the penalized MLE. From now on, we let the event $\cap_{j=1}^4 \mathcal{E}_{n,j}$ hold, where $\mathcal{E}_{n,j}$ for $j = 1, 2, 3, 4$ are defined as in the proof of Theorem 4.1. Recall that these events satisfy

$$\sup_{f \in H^m(1)} P_f^n \left( \cup_{j=1}^4 \mathcal{E}_{n,j}^c \right) \leq (2C_1 + 6) \exp(-Mnhr_n^2).$$

Let $g_n = f - \widehat{f}_{n,\lambda}$. Using Taylor's expansion, we get that

$$PLRT_{n,\lambda}(f) = \ell_{n,\lambda}(f) - \ell_{n,\lambda}(\widehat{f}_{n,\lambda})$$

$$= \ell_{n,\lambda}(\widehat{f}_{n,\lambda} + g_n) - \ell_{n,\lambda}(\widehat{f}_{n,\lambda})$$

$$= \int_0^1 \int_0^1 s DS_{n,\lambda}(\widehat{f}_{n,\lambda} + ss'g_n) g_n g_n ds ds'$$

$$= \int_0^1 \int_0^1 s [DS_{n,\lambda}(\widehat{f}_{n,\lambda} + ss'g_n) - DS_{n,\lambda}(f)] g_n g_n ds ds'$$

$$\quad + \frac{1}{2} [DS_{n,\lambda}(f) - DS_\lambda(f)] g_n g_n + \frac{1}{2} DS_\lambda(f) g_n g_n.$$

For any $s, s' \in \mathbb{I}$, it is easy to see that $\|f - \widehat{f}_{n,\lambda} - ss'g_n\|_{\sup} = (1 - ss')\|g_n\|_{\sup} \leq c_K h^{-1/2}\delta'_n(M) \leq 2/3$ (see (S.16)).

$$
\begin{aligned}
&|[DS_{n,\lambda}(\widehat{f}_{n,\lambda} + ss'g_n) - DS_{n,\lambda}(f)]g_n g_n| \\
\leq\ & |\frac{1}{n}\sum_{i=1}^{n}[\ell^{(2)}(Y_i - \widehat{f}_{n,\lambda} - ss'g_n(X_i)) - \ell^{(2)}(Y_i - f(X_i))]g_n(X_i)^2| \\
\leq\ & \frac{1}{n}\sum_{i=1}^{n}\sup_{|a|\leq 1}|\ell^{(3)}(\eta_i + a)| \cdot \|g_n\|_{\sup}^3 \\
\leq\ & C_0(Mnhr_n^2 + \log n)(c_K h^{-1/2}\delta'_n(M))^3.
\end{aligned}
$$

Define

$$
\psi_{n,f}^{(3)}(Z_i; g) = \frac{\ell^{(2)}(\eta_i)g(X_i)c_K^{-1}h^{1/2}}{C_0(Mnhr_n^2 + \log n)}I_{A_i}, \ g \in \mathcal{G},
$$

where $A_i$ is defined in (i). Then it is easy to see that for any $g_1, g_2 \in \mathcal{G}$,

$$
|\psi_{n,f}^{(3)}(Z_i; g_1) - \psi_{n,f}^{(3)}(Z_i; g_2)| \leq \frac{|\ell^{(2)}(\eta_i)|I_{A_i}}{C_0(Mnhr_n^2 + \log n)}c_K^{-1}h^{1/2}\|g_1 - g_2\|_{\sup} \leq c_K^{-1}h^{1/2}\|g_1 - g_2\|_{\sup}.
$$

It follows by Lemma 6.1 that $\sup_{f\in H^m(1)} P_f^n(\mathcal{E}_{n,5}^c) \leq 2\exp(-Mnhr_n^2)$, where $\mathcal{E}_{n,5} = \{\sup_{g\in\mathcal{G}}\|Z_{n,f}^{(3)}(g)\| \leq \sqrt{Mnhr_n^2}A(h)\}$ and $Z_{n,f}^{(3)}(g) = \frac{1}{\sqrt{n}}\sum_{i=1}^{n}[\psi_{n,f}^{(3)}(Z_i; g)K_{X_i} - E_f\{\psi_{n,f}^{(3)}(Z_i; g)K_{X_i}\}]$. From now on, we assume $\cap_{j=1}^{5}\mathcal{E}_{n,j}$ holds.

It can be seen by (S.13) and independence between $\eta$ and $X$ that

$$
\begin{aligned}
|E_f^Z\{\ell^{(2)}(\eta_i)g_n(X_i)^2 I_{A_i^c}\}| &\leq E\{|\ell^{(2)}(\eta_i)|I_{A_i^c}\}E\{g_n(X)^2\} \\
&\leq \sigma^{-2}E\{|\ell^{(2)}(\eta)|^2\}^{1/2}P_f^n(A_i^c)^{1/2}\|g_n\|^2 \\
&\leq \sigma^{-2}\sqrt{2C_1C_0^2}\sqrt{2C_1}\exp(-(Mnhr_n^2 + \log n)/2)\delta'_n(M)^2 \\
&= 2C_0C_1\sigma^{-2}\exp(-(Mnhr_n^2 + \log n)/2)\delta'_n(M)^2.
\end{aligned}
$$

Note that

$$
\begin{aligned}
&|[DS_{n,\lambda}(f) - DS_\lambda(f)]g_n g_n| \\
=\ & \frac{1}{n}|\sum_{i=1}^{n}[\ell^{(2)}(\eta_i)g_n(X_i)^2 - E_f^Z\{\ell^{(2)}(\eta_i)g_n(X_i)^2\}]| \\
\leq\ & \frac{1}{n}|\sum_{i=1}^{n}[\ell^{(2)}(\eta_i)g_n(X_i)^2 I_{A_i} - E_f^Z\{\ell^{(2)}(\eta_i)g_n(X_i)^2 I_{A_i}\}]| + |E_f^Z\{\ell^{(2)}(\eta_i)g_n(X_i)^2 I_{A_i^c}\}| \\
=\ & |\langle\frac{1}{n}\sum_{i=1}^{n}[\ell^{(2)}(\eta_i)g_n(X_i)I_{A_i}K_{X_i} - E_f^Z\{\ell^{(2)}(\eta_i)g_n(X_i)I_{A_i}K_{X_i}\}], g_n\rangle| + |E_f^Z\{\ell^{(2)}(\eta_i)g_n(X_i)^2 I_{A_i^c}\}| \\
\leq\ & \|\frac{1}{n}\sum_{i=1}^{n}[\ell^{(2)}(\eta_i)g_n(X_i)I_{A_i}K_{X_i} - E_f^Z\{\ell^{(2)}(\eta_i)g_n(X_i)I_{A_i}K_{X_i}\}]\| \cdot \|g_n\| + |E_f^Z\{\ell^{(2)}(\eta_i)g_n(X_i)^2 I_{A_i^c}\}|.
\end{aligned}
$$

Letting $\bar{g} = g_n/d'_n$, where $d'_n = c_K h^{-1/2} \delta'_n(M)$, we get that

$$\|\frac{1}{n}\sum_{i=1}^{n}[\ell^{(2)}(\eta_i)g_n(X_i)I_{A_i}K_{X_i} - E_f^Z\{\ell^{(2)}(\eta_i)g_n(X_i)I_{A_i}K_{X_i}\}]\|$$

$$= \frac{d'_n}{n}\|\sum_{i=1}^{n}[\ell^{(2)}(\eta_i)\bar{g}(X_i)I_{A_i}K_{X_i} - E_f^Z\{\ell^{(2)}(\eta_i)\bar{g}(X_i)I_{A_i}K_{X_i}\}]\|$$

$$= \frac{d'_n}{n}C_0(Mnhr_n^2 + \log n)c_K h^{-1/2}\|\sum_{i=1}^{n}[\psi_{n,f}^{(3)}(Z_i;\bar{g})K_{X_i} - E_f^Z\{\psi_{n,f}^{(3)}(Z_i;\bar{g})K_{X_i}\}]\|$$

$$= C_0(Mnhr_n^2 + \log n)(c_K h^{-1/2})^2 n^{-1/2}\delta'_n(M)\|Z_{n,f}^{(3)}(\bar{g})\|$$

$$\leq \sqrt{M}C_0 c_K^2 (Mnhr_n^2 + \log n)h^{-1/2}\delta'_n(M)r_n A(h),$$

where the last inequality follows by $\bar{g} \in \mathcal{G}$.

By the above analysis, it holds that

$$|[DS_{n,\lambda}(f) - DS_{\lambda}(f)]g_n g_n|$$

$$\leq \sqrt{M}C_0 c_K^2(Mnhr_n^2 + \log n)h^{-1/2}\delta'_n(M)^2 r_n A(h) + 2C_0 C_1 \sigma^{-2}\exp(-(Mnhr_n^2 + \log n)/2)\delta'_n(M)^2.$$

Therefore, for any $f \in H^m(1)$, on $\cap_{j=1}^{5}\mathcal{E}_{n,j}$, it holds that

$$|PLRT_{n,\lambda}(f) + \frac{1}{2}\|g_n\|^2|$$

$$\leq \frac{1}{2}C_0 c_K^3(Mnhr_n^2 + \log n)h^{-3/2}\delta'_n(M)^3 + \frac{1}{2}\sqrt{M}C_0 c_K^2(Mnhr_n^2 + \log n)h^{-1/2}\delta'_n(M)^2 r_n A(h)$$

$$+ C_0 C_1 \sigma^{-2}\delta'_n(M)^2 \exp(-(Mnhr_n^2 + \log n)/2) = R_n.$$

S.7. *Proofs of the auxilliary results.*

S.7.1. *Proof of Lemma S.6.* Since $\epsilon_i^2$ are independent sub-exponential random variables given $(X_1, \ldots, X_n)$. Therefore, by standard tail bounds for sub-exponential random variables [14, Lemma 1], we have that for any $M > 0$

$$(S.19) \qquad P_\epsilon\Big(\Big|\sum_{i=1}^{n}\epsilon_i^2 K(X_i, X_i) - \sum_{i=1}^{n}K(X_i, X_i)\Big| \geq 2g_n\sqrt{M} + 2h^{-1}c_K^2 M\Big) \leq 2\exp(-M),$$

where $g_n^2 = \sum_{i=1}^{n}K^2(X_i, X_i)$. Now we apply Hoeffding's inequality to obtain

$$(S.20) \qquad P_X\Big(\Big|\sum_{i=1}^{n}K(X_i, X_i) - nE\|K_{X_i}\|^2\Big| \geq c_K\sqrt{\frac{nM}{2h}}\Big) \leq 2\exp(-M),$$

$$(S.21) \qquad P_X\Big(\sum_{i=1}^{n}K^2(X_i, X_i) \geq nE\|K_{X_i}\|^4 + \frac{c_K^2}{h}\sqrt{\frac{nM}{2}}\Big) \leq \exp(-M).$$

Since $(X_1, \ldots, X_n)$ and $(\epsilon_1, \ldots, \epsilon_n)$ are independent, we can apply inequalities (S.19), (S.20) and (S.21) to obtain

$$P\Big\{\Big|\sum_{i=1}^n \epsilon_i^2 K(X_i, X_i) - nE\,\|K_{X_i}\|^2\Big| \geq \frac{c_K^2}{h}\sqrt{2nM} + h^{-1/2}c_K\sqrt[4]{2nM^3} + 2h^{-1}c_K^2 M + c_K\sqrt{\frac{nM}{2h}}\Big\}$$
$$\leq 5\exp(-M),$$

where we used the fact that $\|K_{X_i}\| \leq h^{-1/2}c_K$ and the inequality $\sqrt{x+y} \geq \sqrt{x/2} + \sqrt{y/2}$ for $x, y > 0$. This is equivalent to the claimed result since $E\{\|S_{n,\lambda}(0)\|^2\} = n^{-1}E\,\|K_{X_i}\|^2$.

S.7.2. *Proof of Lemma S.7.* By Hanson-Wright inequality [11], given $(X_1, \ldots, X_n)$ for any $M > 0$

$$(S.22) \qquad P_\epsilon\Big(\Big|\sum_{1\leq i<j\leq n} \epsilon_i\epsilon_j K(X_i, X_j)\Big| \geq 2e_n\sqrt{M} + 2f_n M\Big) \leq 2\exp(-M),$$

where $e_n^2 = \sum_{1\leq i<j\leq n} K^2(X_i, X_j)$ and $f_n = \lambda_{\max}(K_n)$. Here for a p.s.d. matrix $A$, $\lambda_{\max}(A)$ stands for its largest eigenvalue, and the kernel matrix $K_n = [K(X_i, X_j)]_{n\times n}$.

By the bounded differences inequality, for any $M > 0$

$$(S.23) \qquad P_X\Big(\sum_{1\leq i<j\leq n} K^2(X_i, X_j) \geq \frac{n^2\rho_K^2}{h} + \frac{4n^2 c_K^2}{h}\sqrt{\frac{M}{2n}}\Big) \leq \exp(-M).$$

By Theorem 4 in [5], for any $x > 0$

$$(S.24) \qquad P_X\Big(\frac{f_n}{n} \geq \mu + 2\sqrt{\frac{E\,\|K_{X_i}\|^4}{n}} + \frac{3c_K^2}{h}\sqrt{\frac{x}{2n}}\Big) \leq 3\exp(-M).$$

Here $\mu$ is the largest eigenvalue of integral operator $K$ on $S^m(\mathbb{I})$ defined as

$$K(u, v) = \sum_{\nu\geq 1} \frac{1}{1 + \lambda\rho_\nu} \phi_\nu(u)\phi_\nu(v)$$

so that

$$Kf(\cdot) = \int_0^1 f(u)K(u, \cdot)dP(u).$$

It is easy to see that $\mu = (1 + \lambda\rho_1)^{-1} \leq 1$, which combined with (S.24) implies

$$(S.25) \qquad P_X\Big(f_n \geq n + 2\sqrt{n\,E\,\|K_{X_i}\|^4} + \frac{3c_K^2}{h}\sqrt{\frac{nx}{2}}\Big) \leq 3\exp(-M).$$

Since $(X_1, \ldots, X_n)$ and $(\epsilon_1, \ldots, \epsilon_n)$ are independent, we can combine inequalities (S.22), (S.24) and (S.25), and use the fact that $E\,\|K_{X_i}\|^4 \leq h^{-2}c_K^4$ to obtain

$$P\Big(\Big|\sum_{1\leq i<j\leq n} \epsilon_i\epsilon_j K(X_i, X_j)\Big| \geq \frac{2n\rho_K}{\sqrt{h}}\sqrt{M} + \frac{4n^{3/4}c_K}{\sqrt[4]{2}\sqrt{h}}M^{3/4} + 2nM + \frac{4\sqrt{n}c_K^2}{h}M + \frac{3c_K^2}{\sqrt{n}}\sqrt{2M}\Big)$$
$$\leq 8\exp(-M),$$

which implies the claimed result.

S.7.3. *Proof of Lemma S.8.* Since $\epsilon_i$ are 1-sub-Gaussian, we have

$$P_\epsilon\Big(\Big|\frac{1}{n}\sum_{i=1}^n \epsilon_i f(X_i)\Big| \geq \|f\|_n\sqrt{\frac{2M}{n}}\Big) \leq 2\exp(-M)$$

given $(X_1,\ldots,X_n)$ for any fixed function $f$ in $S^m(\mathbb{I})$, where $\|f\|_n^2 = n^{-1}\sum_{i=1}^n f^2(X_i)$. By Bernstein's inequality

$$P_X\Big(\|f\|_n^2 \geq \|f\|^2 + \frac{2\|f\|}{n}\sqrt{M} + \frac{2\|f\|_{\sup}}{n}M\Big) \leq \exp(-M).$$

Combining the last two displayed results yields the claimed inequality.

S.7.4. *Proof of Lemma S.9.* Since $f_X^\perp$ is a linear function, we have $J(f_X^\perp, f_X^\perp) = 0$ and

$$
\begin{aligned}
\|f_X^\perp\|^2 &= \|f_X^\perp\|_2^2 \leq 2\|(D_X^T D_X)^{-1}D_X^T f^\perp(X_1^n)\|^2 \leq 2\lambda_X^2\|D_X^T f^\perp(X_1^n)\|^2 \\
&= 2\lambda_X^2\Big\{\Big(\sum_{i=1}^n f^\perp(X_i)\Big)^2 + \Big(\sum_{i=1}^n X_i f^\perp(X_i)\Big)^2\Big\},
\end{aligned}
\tag{S.26}
$$

where $\lambda_X$ is the largest eigenvalue of matrix $(D_X^T D_X)^{-1}$.

According to the proof of Lemma 6.3, $\lambda_X$ has the following tail bound

$$P_X\Big(\lambda_X \geq \frac{2}{n}\Big(1 - \frac{1}{n} - \sqrt{\frac{M}{n}} - \frac{M}{n}\Big)^{-1}\Big) \leq \exp(-M)$$

for all $M \leq n/4$.

Now consider the two terms inside the curly brackets of (S.26). Since $f \in H^m(1)$, we have $\|f^\perp\|_\infty \leq 1$. By definition, we know that $E_X[f^\perp(X)] = E_X[Xf^\perp(X)] = 0$. As a consequence, an application of Hoeffding's inequality yields

$$P_X\Big(\Big|\sum_{i=1}^n f^\perp(X_i)\Big| \geq \sqrt{2n\,M}\Big) \leq 2\exp(-M), \quad \text{and}$$

$$P_X\Big(\Big|\sum_{i=1}^n X_i f^\perp(X_i)\Big| \geq \sqrt{2n\,M}\Big) \leq 2\exp(-M), \quad \text{for all } x > 0.$$

Combining the last three displays and inequality (S.26), we obtain

$$P\Big\{\|f_X^\perp\|^2 \geq \frac{32}{n}M\Big(1 - \frac{1}{n} - \sqrt{\frac{M}{n}} - \frac{M}{n}\Big)^{-2}\Big\} \leq 5\exp(-M),$$

for any $M \leq n/4$.

Department of Statistics
Florida State University
117 N. Woodward Ave.
Tallahassee, FL 32306
Email: yyang@stat.fsu.edu

Department of Mathematical Sciences
Binghamton University
4400 Vestal Parkway East
Binghamton, NY 13902
Email: zshang@binghamton.edu

Department of Statistics
Purdue University
250 N. University Street
West Lafayette, IN 47906
Email: chengg@purdue.edu